\documentclass[leqno,a4paper]{article}

\usepackage{amsmath}
\usepackage{amsfonts}
\usepackage{enumerate}
\usepackage{theorem}

\theoremstyle{plain}
\newtheorem{teo}{Theorem}[section]
\newtheorem{lem}[teo]{Lemma}

\newtheorem{prop}[teo]{Proposition}

{\theorembodyfont{\rmfamily}
\newtheorem{defin}[teo]{Definition}
\newtheorem{oss}[teo]{Remark}
\newtheorem{exam}[teo]{Example}}

\renewcommand{\eqref}[1]{\textnormal{(\ref{#1})}}

\numberwithin{equation}{section}

\newcommand{\cvd}{\hfill$\square$}
\newcommand{\proof}[1]{\noindent\textsc{Proof#1}}
\newcommand{\rmi}{\mathrm{i}}
\newcommand{\rme}{\mathrm{e}}

\title{Examples of exponential instability
for elliptic inverse problems}

\author{Michele Di Cristo\footnotemark[2]\hspace{1em}and\hspace{1em}Luca Rondi\footnotemark[3]\\
\normalsize{Dipartimento di Scienze Matematiche}\\
\normalsize{Universit\`a degli Studi di Trieste, Italy}}

\date{}

\begin{document}

\maketitle
\footnotetext[2]{Supported by MIUR under grant n.~2002013279. E-mail: \texttt{dicristo@mat.unimi.it}}
\footnotetext[3]{Supported by MIUR under grant n.~2002013279 and under Progetto Giovani Ricercatori.
E-mail: \texttt{rondi@mathsun1.univ.trieste.it}}

\setcounter{section}{0}
\setcounter{secnumdepth}{2}

\begin{abstract}
Following a recent paper by N.~Mandache (Inverse Problems \textbf{17} (2001),
pp.~1435--1444),
we establish a general procedure for determining the instability character of
inverse problems. We apply this procedure to many elliptic inverse
problems concerning the
determination of defects of various types
by different kinds of boundary measurements
and we show that these problems are exponentially ill-posed.
\end{abstract}

\section{Introduction}

Many inverse problems associated to partial differential equations concern the problem
of determining a parameter of the equation, for example either a coefficient of the
equation (\emph{coefficient identification})
or the geometry (that is the boundary) of the region where the phenomenon
modelled by the equation occurs (\emph{boundary identification}). In order to determine
this parameter one needs additional information on the solutions to the
partial differential equation, usually constituted of measurements of the solutions
on an accessible (and therefore known) part of the boundary of the region in which the
phenomenon takes place.

As an example of a coefficient identification problem, we may think of the \emph{inverse
conductivity problem}, whose formulation is due to A.~P.~Calder\'on, see \cite{Cal}.
In this problem, electrostatic measurements of voltage and current are collected on
the boundary of a conducting body and
by these data one tries to obtain information about the conductivity inside the body.

For what concerns boundary identification problems, we consider the following examples.
First, determination of a defect inside a conducting body by electrostatic measurements on
the boundary. The defect can be of many different types: it can be an
\emph{inclusion}, that is a region where the conductivity is different from the background
conductivity, see for instance \cite{Ale e Isak} and \cite{Isak88}; it can be a
\emph{crack}, that is a
fracture, as it has been introduced in \cite{Fri e Vog}, see also \cite{Ale e DiB} and
\cite{Ron99:2}; it can be a \emph{cavity} or a
\emph{boundary material loss},
due to corrosion for example, see, for instance,
\cite{Ale e Bere e Ross e Ves} and \cite{Ron99:2}. Then, also the
determination of an \emph{obstacle}
by acoustic measurements in the far-field can be considered as a problem
of this kind, see \cite{Col e Kre98}.

It has been noted several times that these kinds of
inverse problems are ill-posed; in fact, even if the amount of data collected
is sufficient to guarantee uniqueness, the coefficient or the defect, respectively,
usually does
not depend continuously, that is in a stable way, from the measured data.

For the numerical treatment of inverse problems, the ill-posedness constitutes a severe
difficulty. The second main difficulty is
usually due to the fact that inverse problems are tipically non-linear, even if
the direct problem, in the examples above a boundary value problem for
an elliptic partial differential equation, is a linear one.
An accurate knowledge of the character of the ill-posedness is an advantage
for devising efficient numerical methods. Since the problem is ill-posed, in order to
recover some kind of stability, we need to apply a regularization procedure, that is
to restrict the space of admissible unknowns
(either the coefficients or the defects) by assuming that they satisfy a priori conditions
involving usually some kind of smoothness assumptions. With
this a priori information, it is possible to prove that the unknowns depend in a
continuous way from the measured data. However, an explicit knowledge of such a continuous
dependence is crucial for several reasons. First, it provides us with a quantitative
information on how much ill-posed the inverse problem is,
thus how much difficult it is to solve it
numerically; second, a precise knowledge of the modulus of continuity of the dependence
of the unknowns from the measured data indicates the
optimal rate of convergence for regularization schemes and can be useful also for
tuning the regularization parameter, see for instance \cite{Eng e Han e Neu}.

The determination of the modulus of continuity has to be done in two steps. First,
we have to establish \emph{stability estimates}
conditioned to some a priori assumptions on the unknowns; second, we have to show that
these stability estimates are \emph{optimal} or at least
essentially optimal. In order to fulfil
the second part of this program, we need to construct examples
which show that the inverse problem has an instability
character of the same order, or at least of the same kind, that is of logarithmic or
H\"older type, for instance,
of the stability estimates already
established.
We say that our inverse problem is \emph{exponentially ill-posed}, or
\emph{severely ill-posed}, if such a modulus of continuity
is of logarithmic type. In other words, exponential instability corresponds to
the fact that optimal stability estimates are at most of logarithmic type.

The first of these examples has been constructed in \cite{Ale97}
and deals with the problem of the determination of a boundary material loss
in a planar conductor. This example
shows that the stability estimates developed in \cite{Ron99:2} are essentially
optimal; since these estimates are of logarithmic type, this kind of problem is
therefore exponentially ill-posed.
An example similar
to the one in \cite{Ale97} has been constructed for the problem of cavities, still in two
dimensions, in \cite{Ale e Ron99}. These two examples are explicit in the
sense that a family of solution showing the instability character of the problem
is given by explicit formulas,
choosing defects whose boundaries are highly oscillating.
The construction of a family satisfying the instability property looked for
is not an easy task for other inverse problems.

Recently, however, N.~Mandache has proved in \cite{Man} that the
inverse conductivity problem is also exponentially unstable, showing at the same
time that the estimates given in \cite{Ale88:2} are optimal. The procedure used in
\cite{Man} does not depend on an explicit construction, it is instead constituted by a
purely topological argument, which follows from the work of
A.~N.~Kolmogorov and V.~M.~Tihomirov, \cite{Kol e Tih}.
We wish to illustrate the argument as follows. Let
$F:X\mapsto Y$ be a function, $X$ and $Y$ being metric spaces. As a model of an inverse
problem,
$X$ represents the space of unknowns, $Y$ the space of the measured data and $F$ is the
forward map representing the direct problem.
Let us assume that
there exists $x_0\in X$ so that 
for every $\varepsilon>0$ the ball $B(x_0,\varepsilon)$ contains $f(\varepsilon)$
disjoint balls of radius $\varepsilon/2$, $f(\varepsilon)$
being an integer depending on $\varepsilon$. Furthermore, we assume that
for every $\delta>0$ there exists
an integer $g(\delta)$ so that $F(X)$ can be covered by 
$g(\delta)$ balls of radius $\delta$. If we can find $\varepsilon_1>0$ and
$\delta(\varepsilon)$ so that
for every $\varepsilon$, $0<\varepsilon<\varepsilon_1$,
$f(\varepsilon)>g(\delta(\varepsilon))$, then we can find $x_1$ and
$x_2$ in $B(x_0,\varepsilon)$ so that $d_X(x_1,x_2)\geq \varepsilon$ and
$d_Y(F(x_1),F(x_2))\leq 2\delta(\varepsilon)$. Thus, $\delta(\varepsilon)$
provides
an indication of the instability of the inverse to the map $F$. Hence,
it appears clear that establishing this instability
character depends on an accurate counting either of the maximal amount
of disjoint balls with fixed radius that can be found in a given
ball of the space $X$ or of
the minimal amount of balls with fixed radius
required to cover the image through $F$ of $X$.

This procedure immediately appears to be very general and very well suited to be applied
in the context of
inverse problems. In fact, the space of unknowns has, in general,
a richer structure with
respect to that of the data, since usually in inverse problems the
forward map $F$ is compact. A first application of this procedure to ill-posed
problems is developed in \cite{Sca e Via}.

Following the topological arguments of \cite{Kol e Tih} and the procedure described in
\cite{Man}, we have extracted a general method for determining instability,
to be applicable to many different inverse problems. In 
Theorem~\ref{abstractthm} below, we have stated in a rather abstract
framework the outline of the method, in one of its possible formulations (for slightly
different but analogous formulations we refer, for instance, to the discussion of
the inverse scattering case, see Subsection~\ref{scatsubproof}).

Then we have applied our abstract result to many inverse boundary value problems
of elliptic type. We have shown that all the kinds of boundary identification problems
briefly described above are exponentially ill-posed.
We also wish to remark that, as in the explicit examples of
\cite{Ale97} and \cite{Ale e Ron99}, the ill-posedness is of exponential type no matter
which and how many measurements we take. Our examples in fact deal with the ideal case
of performing all possible measurements. This fact is somewhat surprising since
in these boundary identification problems
a much lesser amount of data is required to have unique identification of a defect
and also to have stability estimates; usually a finite number of measurements is enough.
This shows the difficulty of the problem and that performing more measurements
or different ones does not solve the problem of ill-posedness.

The plan of the paper is as follows. In order to point out to the reader all the inverse problems
to which we have successfully applied the method,
we first describe, in Section~\ref{ressec},
all the instability results that are contained in the paper. Then we proceed with the
proofs of these results. The proofs are divided into three sections. 
In Section~\ref{abstractsec} we state and prove an abstract result, Theorem~\ref{abstractthm},
which provides the general procedure
for obtaining the instability examples and therefore constitutes the key ingredient and crucial part of the proofs
of all the results
described in Section~\ref{ressec}. In fact, the proofs are in general obtained as straightforward applications
of this abstract theorem. In order to apply the abstract theorem, what is essentially needed is to choose a suitable
orthonormal basis and to check that all the hypotheses of the abstract theorem are satisfied. Concerning
orthonormal basis, we
shall employ eigenfunctions corresponding to eigenvalue problems of Stekloff type. We have collected all the information
we shall need about these orthonormal basis in Section~\ref{prelsec}. 
Then, in Section~\ref{examsecproof}, the proofs of the instability results are concluded. Using the orthonormal basis
introduced in Section~\ref{prelsec}, we verify that the abstract result applies to the problems we consider and we
prove their exponential instability.

In details, in Section~\ref{ressec}, first we need to introduce some notations
which will be used repeatedly in the paper. In particular we define, and
investigate the structure of, the metric spaces of the unknowns.
Then, we list the problems for which we have obtained the instability examples, together with
the precise formulation of the instability results. We observe that, for the sake of brevity,
we usually refer to the bibliography
for a more detailed description of the problems considered.
We begin with the problem
of determination of defects of different types by electrostatic boundary measurements.
In Subsection~\ref{inclsub} we treat the problem
of determination of an inclusion, in
Subsection~\ref{cracksub} the determination of cracks is considered,
in Subsection~\ref{cavssub} we deal with the inverse problem of cavities,
in Subsection~\ref{surfcracksub} we treat the case of cracks reaching the boundary
of the domain, that is surface cracks, and in Subsection~\ref{corrsub}
we study the problem of a boundary material loss. Finally, in Subsection~\ref{scatsub},
we deal with inverse scattering
problems, in particular with the determination of obstacles (either of sound-soft or
of sound-hard type) by far-field acoustic measurements. In Section~\ref{abstractsec}
the abstract result is stated and proved.
In Section~\ref{prelsec} we study two different
eigenvalue problems of Stekloff type and we investigate the asymptotic
properties of either their eigenvalues or eigenfunctions, in particular this is done
for three different domains of our interest where the solutions can be computed almost
explicitly.
In Section~\ref{examsecproof}, the conclusions of proofs of all the instability results are developed.

 \section{Statement of the instability results}\label{ressec}

Before stating the main results, we need to introduce some notations about
the Sobolev spaces we shall use and to describe the spaces of the unknowns.

For any $N\geq 2$, any $x=(x_1,\ldots,x_N)\in\mathbb{R}^N$ and any $r>0$, we denote
$B_N(x,r)=\{y\in\mathbb{R}^N:\ \|y-x\|<r\}$. We set
$S^{N-1}(x,r)=\partial B_N(x,r)$. Furthermore, we set
$S^{N-1}=\partial B_N(0,1)$,
and $S^{N-1}_+=\{y\in S^{N-1}:\ y_N\geq 0\}$, and, analogously,
$S^{N-1}_-=\{y\in S^{N-1}:\ y_N\leq 0\}$. Finally, we denote
$B'_{N-1}(x,r)=\{y\in B_N(x,r):\ y_N=x_N\}$.

We need, furthermore, to introduce the following definition.

\begin{defin}\label{def-net-discrete}
Let $(Y,d_Y)$ be a metric space. For a given positive $\delta$, $Y_1$, a subset of
$Y$, is said to be a
$\delta$-\emph{net} for $Y$ if for every $y\in Y$ there exists $y_1\in Y_1$
so that $d_Y(y,y_1)\leq\delta$.

Given $\varepsilon$ positive, $Y_2\subset Y$ is $\varepsilon$-\emph{discrete}
if for any two distinct points $y_2$, $y'_2$ in $Y_2$ we have
$d_Y(y_2,y'_2)\geq\varepsilon$.
\end{defin}

\subsubsection{Notations on Sobolev spaces}\label{sobnot}

Let $\Omega\subset\mathbb{R}^N$, $N\geq 2$, be a bounded domain and let
$\partial\Omega$ be its boundary. About regularity, we assume that
there exists a homeomorphism $\chi:B_N(0,1)\mapsto \Omega$ such that,
for a positive constant $C$, we have
\begin{equation}\label{Lipdomain}
\begin{array}{ll}
\|\chi(\tilde{x})-\chi(\tilde{y})\|\leq C\|\tilde{x}-\tilde{y}\|
& \text{for any }\tilde{x},\tilde{y}\in B_N(0,1),\\
\|\chi^{-1}(x)-\chi^{-1}(y)\|\leq C\|x-y\| & \text{for any }x,y\in\Omega.
\end{array}
\end{equation}

Furthermore, we shall consider two internally disjoint subsets of
$\partial\Omega$, $\Gamma_A$
and $\Gamma_I$, so that $\Gamma_A\cup\Gamma_I=\partial\Omega$.
We assume either that $\Gamma_A=\partial\Omega$ and $\Gamma_I=\emptyset$, or 
that $\Gamma_A$ and $\Gamma_I$ are not empty and 
are assumed to be regular enough,
namely there exists a homeomorphism $\chi:B_N(0,1)\mapsto \Omega$ satisfying
\eqref{Lipdomain}, so that, if we still denote with $\chi$ its extension by continuity
to $\overline{B_N(0,1)}$, then $\Gamma_A=\chi(S^{N-1}_+)$ and
$\Gamma_I=\chi(S^{N-1}_-)$.

We introduce the following
Sobolev spaces. Let $H^1(\Omega)=\{u\in L^2(\Omega):\ \nabla u\in L^2(\Omega)\}$,
where $\nabla u$ denotes the gradient of $u$ in the sense of distributions. We recall
that $H^1(\Omega)$ is a Hilbert space with scalar product
$(u,v)_{H^1(\Omega)}=\int_\Omega\nabla u\cdot\nabla v+uv$.
With $H^{1/2}(\Gamma_A)$ we denote the space of traces of $H^1(\Omega)$ functions
on $\Gamma_A$, which can be endowed in a canonical way
with a scalar product induced by the one of
$H^1(\Omega)$ so that $H^{1/2}(\Gamma_A)$ is a Hilbert space.
By $H^{-1/2}(\Gamma_A)$ we shall denote
the dual space to $H^{1/2}(\Gamma_A)$.
We recall that $H^{1/2}(\Gamma_A)\subset L^2(\Gamma_A)\subset
H^{-1/2}(\Gamma_A)$.
We shall also make use of the following spaces. Let
${}_0H^{1/2}(\Gamma_A)=\{\psi\in H^{1/2}(\Gamma_A):\ \int_{\Gamma_A}\psi=0\}$.
Its dual is given by the space
${}_0H^{-1/2}(\Gamma_A)=\{\eta\in H^{-1/2}(\Gamma_A):\ \langle\eta,1\rangle=0\}$,
where $\langle\cdot,\cdot\rangle$ denotes the duality pairing.

If $\Gamma_I$
is not empty, we set
$H^1_0(\Omega,\Gamma_I)$ and $H^1_{const}(\Omega,\Gamma_I)$
as the closed subspaces of $H^1(\Omega)$ constituted by
the functions $u\in H^1(\Omega)$ so that $u=0$ in a weak sense on $\Gamma_I$
and $u=constant$ in a weak sense on $\Gamma_I$, respectively.
With
$H^{1/2}_0(\Gamma_A,\Omega)$ and $H^{1/2}_{const}(\Gamma_A,\Omega)$
we denote the closed subspaces of
$H^{1/2}(\Gamma_A)$ constituted by the traces of $H^1_0(\Omega,\Gamma_I)$
and $H^1_{const}(\Omega,\Gamma_I)$ functions
on $\Gamma_A$, respectively.

For our purposes, we need to introduce on the Sobolev spaces defined above suitable
scalar products, which are different but topologically equivalent to the
canonical ones. We wish to remark that the definitions of these scalar products 
do not take into account the fact that the spaces $H^{-1/2}$ and $H^{1/2}$
are dual one to each other.

For any $\psi$, $\varphi\in H^{1/2}(\Gamma_A)$, we set $\tilde{\psi}\in H^1(\Omega)$
as the solution to
\begin{equation}
\left\{\begin{array}{ll}
\Delta \tilde{\psi}=0 & \text{in }\Omega,\\
\tilde{\psi}=\psi & \text{on }\Gamma_A,\\
\frac{\partial \tilde{\psi}}{\partial \nu}=0 & \text{on }\Gamma_I,
\end{array}\right.
\end{equation}
and $\tilde{\varphi}$ as the solution to the same boundary value problem with
$\psi$ replaced by $\varphi$, and the scalar product we use on $H^{1/2}(\Gamma_A)$
is given by
\begin{equation}\label{scpr1/2}
(\psi,\varphi)_{H^{1/2}(\Gamma_A)}=
\int_{\Omega}\nabla\tilde{\psi}\cdot\nabla\tilde{\varphi}+\int_{\Gamma_A}\psi\varphi.
\end{equation}
We observe that ${}_0H^{1/2}(\Gamma_A)$ coincides with the subspace which is
orthogonal, with respect to this scalar product,
to the constant function $1$.

Any $\eta\in H^{-1/2}(\Gamma_A)$ can be decomposed, in a unique way, into
the sum of $\hat{\eta}$, an element of ${}_0H^{-1/2}(\Gamma_A)$,
and a constant function $c(\eta)$. Furthermore, to $\hat{\eta}$
we can associate $\tilde{\eta}\in H^1(\Omega)$ that solves
\begin{equation}
\left\{\begin{array}{ll}
\Delta \tilde{\eta}=0 & \text{in }\Omega,\\
\frac{\partial \tilde{\eta}}{\partial \nu}=\hat{\eta} & \text{on }\Gamma_A,\\
\frac{\partial \tilde{\eta}}{\partial \nu}=0 & \text{on }\Gamma_I.
\end{array}\right.
\end{equation}
If, in the same way, we associate to 
$\phi\in H^{-1/2}(\Gamma_A)$ the functions
$\hat{\phi}$, $c(\phi)$ and $\tilde{\phi}$, then the scalar product on
$H^{-1/2}(\Gamma_A)$ may be defined as
\begin{equation}\label{scpr-1/2}
(\eta,\phi)_{H^{-1/2}(\Gamma_A)}=
\int_{\Omega}\nabla\tilde{\eta}\cdot\nabla\tilde{\phi}+c(\eta)c(\phi).
\end{equation}
We remark that, with respect to this scalar product,
${}_0H^{-1/2}(\Gamma_A)$ is the orthogonal subspace to the constant function $1$.

We take $\Gamma_I$ not empty.
If $\psi$ belongs to $H^{1/2}_{const}(\Gamma_A,\Omega)$, then there exist (and are unique)
$\hat{\psi}\in H^{1/2}_0(\Gamma_A,\Omega)$ and a constant function $c(\psi)$ so that
$\psi=\hat{\psi}+c(\psi)$. Let
$\tilde{\psi}\in H^1(\Omega)$
solve
\begin{equation}
\left\{\begin{array}{ll}
\Delta \tilde{\psi}=0 & \text{in }\Omega,\\
\tilde{\psi}=\hat{\psi} & \text{on }\Gamma_A,\\
\tilde{\psi}=0 & \text{on }\Gamma_I.
\end{array}\right.
\end{equation}
Then, if we associate to $\varphi\in H^{1/2}_{const}(\Gamma_A,\Omega)$
its corresponding decomposition given by $\hat{\varphi}$ and $c(\varphi)$,
and its corresponding function $\tilde{\varphi}$,
on $H^{1/2}_{const}(\Gamma_A,\Omega)$ we introduce the scalar product
\begin{equation}\label{scpr1/2const}
(\psi,\varphi)_{H^{1/2}_{const}(\Gamma_A,\Omega)}=
\int_{\Omega}\nabla\tilde{\psi}\cdot\nabla\tilde{\varphi}
+c(\psi)c(\varphi).
\end{equation}
Such a scalar product obviously induces a scalar product on
$H^{1/2}_0(\Gamma_A,\Omega)$, which is the closed subspace of
$H^{1/2}_{const}(\Gamma_A,\Omega)$ orthogonal to the constant function $1$.

\subsubsection{Spaces of smooth perturbations of a given set}

We shall consider the following examples.
Let us fix integers $N\geq 2$ and
$m\geq 1$ and positive constants $\varepsilon$ and $\beta$.
Let us also fix $x\in\mathbb{R}^N$ and $r>0$.

To any real function $f$ defined on $\overline{B'_{N-1}(x,r)}$,
where $B'_{N-1}(x,r)=\{y\in B_N(x,r):\ y_N=x_N\}$,
we associate its \emph{graph}, that is
$\mathrm{graph}(f)=
\{y\in\mathbb{R}^N:\ y_N=
f(y_1,\ldots,y_{N-1},x_N),\ (y_1,\ldots,y_{N-1},x_N)\in\overline{B'_{N-1}(x,r)}\}$,
and,
assuming $f\geq x_N$, its \emph{subgraph}, that is
$\mathrm{subgraph}(f)=
\{y\in\mathbb{R}^N:\ x_N\leq y_N\leq f(y_1,\ldots,y_{N-1},x_N),\
(y_1,\ldots,y_{N-1},x_N)\in\overline{B'_{N-1}(x,r)}\}$.

With the notation $X_{m\beta\varepsilon}(B'_{N-1}(x,r))$ we indicate the space
$\{\mathrm{graph}(f):\ f\in C^m_0(B'_{N-1}(x,r)),
\ \|f\|_{C^m(B'_{N-1}(x,r))}\leq \beta
\text{ and }x_N\leq f\leq x_N+\varepsilon\}$ and with
$Y_{m\beta\varepsilon}(B'_{N-1}(x,r))$ we indicate the space
obtained by taking the subgraphs of all the functions belonging to the same class as
before.
We consider the spaces
$X_{m\beta\varepsilon}(B'_{N-1}(x,r))$ and
$Y_{m\beta\varepsilon}(B'_{N-1}(x,r))$ as metric spaces with the Hausdorff distance.

To any strictly positive function $g$ defined on
$S^{N-1}(x,r)=\partial B_N(x,r)$, we denote its \emph{radial graph} as
$\mathrm{graph}_{rad}(g)=
\{y\in\mathbb{R}^N:\ y=x+g(\omega)\cdot\left(\frac{\omega-x}{r}\right),\ \omega\in
S^{N-1}(x,r)\}$ and its \emph{radial subgraph} as
$\mathrm{subgraph}_{rad}(g)=\{y\in\mathbb{R}^N:\
y=x+\rho\cdot\left(\frac{\omega-x}{r}\right),\ 0\leq \rho\leq g(\omega),\ \omega\in
S^{N-1}(x,r)\}$.

Then, with the notation
$X_{m\beta\varepsilon}(S^{N-1}(x,r))$ we denote the space given by
$\{\mathrm{graph}_{rad}(g):\ g\in C^m(S^{N-1}(x,r)),
\ \|g\|_{C^m(S^{N-1}(x,r))}\leq \beta
\text{ and }r\leq g\leq r+\varepsilon\}$ and with
$Y_{m\beta\varepsilon}(S^{N-1}(x,r))$ we denote the space
of radial subgraphs of all the functions belonging to the same class used
before.
Also the spaces $X_{m\beta\varepsilon}(S^{N-1}(x,r))$
and $Y_{m\beta\varepsilon}(S^{N-1}(x,r))$ are metric spaces endowed
with the Hausdorff distance.

It is an easy remark the fact
that $X_{m\beta\varepsilon}(B'_{N-1}(x,r))$
and $X_{m\beta\varepsilon}(S^{N-1}(x,r))$ are contained in the closed
ball, with respect to the Hausdorff distance between closed sets,
of radius $\varepsilon$
centred at $\overline{B'_{N-1}(x,r)}$ and
$S^{N-1}(x,r)$, respectively. Analogously, $Y_{m\beta\varepsilon}(B'_{N-1}(x,r))$
and $Y_{m\beta\varepsilon}(S^{N-1}(x,r))$ are contained in the closed
ball, again with respect to the Hausdorff distance,
of radius $\varepsilon$ and centre
$\overline{B'_{N-1}(x,r)}$ and
$\overline{B_N(x,r)}$,
respectively.
Maybe more interesting and significant is the fact that the elements of 
$Y_{m\beta\varepsilon}(S^{N-1}(x,r))$ are all compact subsets which are star-shaped with
respect to a common point $x\in\mathbb{R}^N$. The determination of star-shaped sets is
usually considered to be more stable than the determination of other kinds of sets. Nevertheless
many of our examples show that even with a star-shapedness assumption
the instability is still of exponential type.

We would like to study properties of $\varepsilon$-discrete sets of 
$X_{m\beta\varepsilon}(B'_{N-1}(x,r))$, $Y_{m\beta\varepsilon}(B'_{N-1}(x,r))$
and
$X_{m\beta\varepsilon}(S^{N-1}(x,r))$,
$Y_{m\beta\varepsilon}(S^{N-1}(x,r))$. We have the
following proposition.

\begin{prop}\label{discretesetprop}
Let us fix integers $N\geq 2$ and
$m\geq 1$ and positive constants $\beta$ and $r$.
We also fix $x\in\mathbb{R}^N$.
Fixed $\varepsilon>0$, let $X_{\varepsilon}$ be equal to one of the following
four metric spaces\textnormal{:}
$X_{m\beta\varepsilon}(B'_{N-1}(x,r))$, $Y_{m\beta\varepsilon}(B'_{N-1}(x,r))$,
$X_{m\beta\varepsilon}(S^{N-1}(x,r))$ or
$Y_{m\beta\varepsilon}(S^{N-1}(x,r))$.

Then, there exists a positive constant $\varepsilon_0$, depending on
$N$, $m$, $\beta$ and $r$
only, so that for any
$\varepsilon$, $0<\varepsilon<\varepsilon_0$, we can find
$Z_{\varepsilon}$
satisfying the
following properties.
We have that 
the set
$Z_{\varepsilon}$ is
contained in $X_{\varepsilon}$\textnormal{;} $Z_{\varepsilon}$ is
$\varepsilon$-discrete, with respect to the Hausdorff distance\textnormal{;}
and, finally,
$Z_{\varepsilon}$ has at
least $\exp(2^{-N}\varepsilon_0^{(N-1)/m}\varepsilon^{-(N-1)/m})$ elements.
\end{prop}

\proof. The proof can be obtained, with slight modifications, along the lines of
the proof of Lemma~2 in \cite{Man}.\cvd

\subsection{Inverse inclusion problem}\label{inclsub}

Let us assume that the domain $\Omega=B_N(0,1)$, $N\geq 2$,
is occupied by a conducting body.
Let us assume that an \emph{inclusion} $D$ is present
inside the otherwise homogeneous conductor;
that is, there exist two different positive constants $a$ and $b$ 
and a set $D$ which is compactly contained in $\Omega$ (that is $\overline{D}$ is a compact
subset of $\Omega$) so that
the conductivity inside $D$ is constantly equal to $a$ and the conductivity
outside $D$, that is in $\Omega\backslash D$, is constantly equal to $b$.
For the sake of simplicity, we normalize $b$ so that $b=1$ and we take $a$ to be positive
and different from $1$.

The electrostatic potential $u$ inside $\Omega$ is a solution to the following partial
differential equation
\begin{equation}\label{incleq}
\mathrm{div}((1+(a-1)\chi_D)\nabla u)=0\quad\text{in }\Omega,
\end{equation}
where $\chi_D$ denotes the characteristic function of the domain $D$.

Furthermore, $u$ satisfies a boundary condition on $\partial\Omega=S^{N-1}$
which depends on whether we prescribe the voltage $\psi\in H^{1/2}(\partial\Omega)$
on the boundary or we assign the current density
$\eta\in{}_0H^{-1/2}(\partial\Omega)$ on the boundary.
Namely, in the first case the boundary condition is given by
\begin{equation}\label{Dircond}
u=\psi \text{ on }\partial\Omega;
\end{equation}
in the second case by
\begin{equation}\label{Neumcond}
\frac{\partial u}{\partial\nu}=\eta\text{ on }\partial\Omega;
\quad\int_{\partial\Omega}u=0;
\end{equation}
where we have added a normalization condition.

We have existence and uniqueness of a (weak) solution for
both the boundary value problems \eqref{incleq}-\eqref{Dircond} and
\eqref{incleq}-\eqref{Neumcond}.

The inverse problem we consider is the one of recovering the shape and the location
of an unknown inclusion $D$, by performing current and voltage measurements at the
boundary, that is
either by prescribing voltages and measuring the corresponding current densities or
viceversa.

In the literature, a lot of attention has been devoted to the determination of
$D$ by a single measurement; in this case the problem has been often referred to
as the inverse conductivity problem with one measurement. A global uniqueness result
is still missing, see \cite{Ale e Isak} and its references for a more detailed
discussion on this topic.  We remark that if all possible measurements are performed, then the inclusion
can be uniquely determined, see \cite{Isak88}. However,
up to our knowledge, even if all measurements are considered, no explicit stability
estimate for this problem has been established.

We produce an example showing
that, even if we make
many measurements, actually all possible
measurements, the optimal stability for this inverse problem is at the best of
logarithmic type. We treat the case when voltages are prescribed and currents are
measured
and the case in which current densities are assigned and voltages are measured, as well.

We fix two positive integers, $m$ and $N$, $N\geq 2$,
and two positive constants, $\beta$ and $a$, $a\neq 1$. We consider
the metric space $(X,d)$ where $X=Y_ {m\beta(1/4)}(S^{N-1}(0,1/2))$ and
$d$ is the Hausdorff distance. If $D$ belongs to $X$, then we can define
the following two operators.

The operator $\Lambda(D):H^{1/2}(\partial\Omega)\mapsto H^{-1/2}(\partial\Omega)$
is defined as
$$\langle \Lambda(D)\psi,\varphi\rangle=
\langle\frac{\partial u}{\partial\nu}|_{\partial\Omega},\varphi\rangle=
\int_{\Omega}(1+(a-1)\chi_D)\nabla u\cdot\nabla\tilde{\varphi}$$
where $\psi$, $\varphi\in H^{1/2}(\partial\Omega)$, $u$ is the solution to
\eqref{incleq}-\eqref{Dircond} and $\tilde{\varphi}$ is any
$H^1(\Omega)$ function whose trace on $\partial\Omega$ is equal to $\varphi$.
Since the operator $\Lambda(D)$ associates the Dirichlet datum to the
corresponding Neumann datum,
it is usually called the \emph{Dirichlet-to-Neumann} map.

Viceversa, the operator $\mathcal{N}(D):{}_0H^{-1/2}(\partial\Omega)
\mapsto{}_0H^{1/2}(\partial\Omega)$
is given by
$$\mathcal{N}(D)\eta=u|_{\partial\Omega}$$
where $\eta\in{}_0H^{-1/2}(\partial\Omega)$ and $u$ is the solution to
\eqref{incleq}-\eqref{Neumcond}. For the same reasons, the map $\mathcal{N}(D)$ is called
the \emph{Neumann-to-Dirichlet} map.

It is easy to show that for any $D\in X$,
the maps $\Lambda(D)$ and $\mathcal{N}(D)$ are linear and bounded operators
between a Hilbert space and its dual. In the sequel, their norms will be always
assumed to be the canonical ones as bounded operators between Hilbert spaces.
We state the instability result.

\begin{prop}\label{inclprop}
We fix integers $N\geq 2$ and
$m\geq 1$ and a positive constant $\beta$. We also fix $0<a\neq 1$.  
Let $(X,d)$ be a metric space, with
$X=Y_{m\beta(1/4)}(S^{N-1}(0,1/2))$ and $d$ being the Hausdorff distance.
Then we can find positive constants $\varepsilon_1$ and $C$, which
depend on $N$, $m$, $\beta$ and $a$ only,
so that for any $\varepsilon$, $0<\varepsilon<\varepsilon_1$,
there exist $D_1$ and $D_2$ belonging to $X$
such that
\begin{equation}\label{firstpart}
d(D_j,\overline{B_{N}(0,1/2)})\leq\varepsilon,\text{ for any }j=1,2;\quad
d(D_1,D_2)\geq\varepsilon;
\end{equation}
and
\begin{equation}\label{DNpart}
\|\Lambda(D_1)-\Lambda(D_2)\|\leq C\exp(-\varepsilon^{-(N-1)/(2mN)});\\
\end{equation}
\begin{equation}\label{NDpart}
\|\mathcal{N}(D_1)-\mathcal{N}(D_2)\|\leq C\exp(-\varepsilon^{-(N-1)/(2mN)}).
\end{equation}
\end{prop}

The proof of this proposition is postponed to Subsection~\ref{inclsubproof}.

\subsubsection{Experimental measurements}

In Proposition~\ref{inclprop},
the inverse inclusion problem is stated to be
exponentially ill-posed even if we perform all possible measurements of current and voltage type
at the boundary. It is not surprising, therefore, that the inverse inclusion problem is exponentially ill-posed
also with respect to measurements which can be actually obtained from the experiments. We shall refer to this
kind of measurements as the \emph{experimental measurements}. The model which we shall follow
is the one developed in \cite{Som e Che e Isa}, which we briefly describe, referring to the original paper
for more details.

The model is the following. On the boundary of the conductor $\Omega$, we attach $L$ \emph{electrodes}.
The \emph{contact regions} between the electrodes and the conductor are subsets of
$\partial \Omega$ and will
be denoted
by $e_l$, $l=1,\ldots,L$. We assume that the subsets $e_l$, $l=1,\ldots,L$,
are open, connected, with a smooth boundary and so that
their closures are pairwise disjoint.
We remark that we identify any electrode with its contact region.
A current is sent to the body through the electrodes and the corresponding voltages
are measured on the same electrodes. For each $l$, $l=1,\ldots,L$,
the current applied to the electrode $e_l$ will be denoted by
$I_l$ and the voltage measured on the electrode will be denoted by $V_l$. The column vector $I$ whose
components are $I_l$, $l=1,\ldots,L$, is a \emph{current pattern}
if the condition $\sum_{l=1}^LI_l=0$ is satisfied.
The corresponding \emph{voltage pattern}, that is the column vector $V$ whose components are $V_l$, $l=1,\ldots,L$,
is determined up to an additive constant and we always choose to normalize it in such a way that
$\sum_{l=1}^LV_l=0$. The voltage pattern depends on the current pattern in a linear way, through an
$L\times L$ symmetric matrix $R$ which is called the \emph{resistance matrix}, that is
$V=RI$.

The following model
can be used to determine the resistance matrix $R$. We assume that at each
electrode $e_l$, $l=1,\ldots,L$, a \emph{surface impedance} is present and we denote it with $z_l$.
Let us assume that there exists $Z>0$ so that
for each $l$, $l=1,\ldots,L$, $z_l\geq Z$.
Let $D$ be as before an inclusion in $\Omega$. The conductivity in $D$ is $a$, where $a$ is a
positive constant different from $1$, and the conductivity outside $D$ is equal to $1$. 
If we apply the current pattern
$I$ on the electrodes, then the voltage $u$ inside the body satisfies the following
boundary value problem
\begin{equation}\label{expmeaspbm}
\left\{\begin{array}{ll}
\mathrm{div}((1+(a-1)\chi_D)\nabla u)=0 &\text{in }\Omega,\\
u+z_l\frac{\partial u}{\partial \nu}=U_l &\text{on }e_l,\ l=1,\ldots,L,\\
\frac{\partial u}{\partial \nu}=0 &\text{on }\partial\Omega\backslash\bigcup_{l=1}^Le_l,\\
\int_{e_l}\frac{\partial u}{\partial \nu}=I_l&\text{for any }l=1,\ldots,L,
\end{array}\right.
\end{equation}
where $U_l$, $l=1,\ldots,L,$ are constants to be determined. We call $U$ the column vector whose
components are given by $U_l$, $l=1,\ldots,L$.

For any $l$, $l=1,\ldots,L$, $V_l$, a component of the voltage pattern $V$, is given by
$V_l=\int_{e_l}u$, thus, by \eqref{expmeaspbm},
$$V_l=|e_l|U_l-z_lI_l,$$
where $|e_l|$ denotes the surface measure of $e_l$.

By \cite[Theorem~3.3]{Som e Che e Isa}, we infer that there exists a unique couple $(u,U)$,
$u$ being in $H^1(\Omega)$ and $U$ being a column vector with $L$ components so that
$\sum_{l=1}^L|e_l|U_l-z_lI_l=0$, such that \eqref{expmeaspbm} is satisfied. Thus the current pattern
$I$ uniquely determines the voltage pattern $V$, if this is normalized in such a way that
$\sum_{l=1}^LV_l=0$. Furthermore, it has been proved in \cite{Som e Che e Isa} that the relation
between $I$ and $V$ is linear, thus the resistance matrix $R(D)$ is well defined. Finally, it has
been shown that $R(D)$ is actually symmetric. We remark that we shall assume, without loss
of generality, that $R(D)[1]=0$, where $[1]$ denotes the column vector whose components are all equal to $1$.
Also, we recall that the norm of $R(D)$ will always be the norm of linear operators from $\mathbb{R}^L$ into itself.

The following instability result will be proved in Subsection~\ref{inclsubproof}.

\begin{prop}\label{inclpropexpmeas}
We fix integers $N\geq 2$ and
$m\geq 1$ and a positive constant $\beta$. We also fix $0<a\neq 1$ and $Z>0$.  
Let $(X,d)$ be a metric space, with
$X=Y_{m\beta(1/4)}(S^{N-1}(0,1/2))$ and $d$ being the Hausdorff distance.
Let us assume that $L\geq 2$ electrodes $e_l$,
$l=1,\ldots,L$, and their surface impedances $z_l$, $l=1,\ldots,L$, are fixed
and satisfy the previously stated assumptions.
Then we can find positive constants $\varepsilon_1$ and $\tilde{C}$, which
depend on $N$, $m$, $\beta$, $a$, $Z$ and the electrodes only,
so that for any $\varepsilon$, $0<\varepsilon<\varepsilon_1$,
there exist $D_1$ and $D_2$ belonging to $X$
such that
\begin{equation}
\begin{array}{l}
d(D_j,\overline{B_{N}(0,1/2)})\leq\varepsilon,\text{ for any }j=1,2;\quad
d(D_1,D_2)\geq\varepsilon;\\
\|R(D_1)-R(D_2)\|\leq \tilde{C}\exp(-\varepsilon^{-(N-1)/(2mN)}).
\end{array}
\end{equation}
\end{prop}

\subsection{Inverse crack problem}\label{cracksub}

Let $\Omega=B_N(0,1)$, $N\geq 2$, be the region occupied by a
homogeneous conducting body.
Let us assume that inside the conductor there is a \emph{crack} $\sigma$, that is a
closed set inside $\Omega$ so that $\Omega\backslash\sigma$ is connected and, locally,
$\sigma$
can be represented by the graph of a smooth function.
We can consider two different types of cracks, \emph{perfectly insulating} and
\emph{perfectly conducting}, and we can prescribe on the (exterior) boundary of $\Omega$
either the voltage or the current density. Thus, the electrostatic potential $u$ in
$\Omega$ satisfies either
\begin{equation}\label{inscrack}
\left\{
\begin{array}{ll}
\Delta u=0 & \text{in }\Omega\backslash\sigma,\\
\frac{\partial u}{\partial \nu}=0 & \text{on }\partial\sigma;
\end{array}
\right.
\end{equation}
if $\sigma$ is perfectly insulating, or, when $\sigma$ is assumed to be perfectly
conducting,
\begin{equation}\label{condcrack}
\left\{
\begin{array}{ll}
\Delta u=0 & \text{in }\Omega\backslash\sigma,\\
u=constant & \text{on }\sigma.
\end{array}
\right.
\end{equation}
We remark that, in \eqref{inscrack}, on $\partial\sigma$ means on either sides
of $\sigma$.
On the boundary the potential satisfies either
\begin{equation}\label{voltbdrycond}
u=\psi \text{ on }\partial\Omega;\quad
\langle\frac{\partial u}{\partial \nu}|_{\partial\Omega},1\rangle=0;
\end{equation}
where $\psi\in H^{1/2}(\partial\Omega)$ is the prescribed voltage at the boundary,
or, if we prescribe the current density on the boundary to be
$\eta\in{}_0H^{-1/2}(\partial\Omega)$,
\begin{equation}\label{currentbdrycond}
\frac{\partial u}{\partial \nu}=\eta \text{ on }\partial\Omega;\quad
\int_{\partial\Omega}u=0;
\end{equation}
we wish to remark that normalization conditions have been added to the
boundary conditions.

We have that all the direct problems \eqref{inscrack}-\eqref{voltbdrycond},
\eqref{inscrack}-\eqref{currentbdrycond}, \eqref{condcrack}-\eqref{voltbdrycond}
and \eqref{condcrack}-\eqref{currentbdrycond} admit a unique (weak) solution.

The inverse crack problem consists of recovering the shape and location of an unknown
crack $\sigma$ by performing electrostatic measurements at the boundary.

In this subsection we shall state the instability character of such 
an inverse problem, in all the possible cases, that is when we consider either
insulating
or conducting cracks, and when either we prescribe voltages and measure
corresponding currents or
we prescribe currents and measure corresponding voltages.

For a detailed analysis of uniqueness and stability of this problem we refer to
\cite{Ron99:2}, for the two-dimensional case, and to \cite{Ale e DiB}, for
the three-dimensional case, and to their bibliographies. We wish to remark that,
for what concerns uniqueness and stability results, these have been obtained
with a finite number of boundary measurements, usually with two suitably chosen
measurements. Our instability example shows the optimality of the stability
estimates previously obtained and that the stability can not be improved
by taking different or more measurements.

The framework of our
example is as follows. Let $N\geq 2$ and $m$, positive integers,
and $\beta$, a positive constant, be fixed.
Let $X=X_{m\beta(1/4)}(B'_{N-1}(0,1/2))$ with the Hausdorff distance.
To any $\sigma\in X$, we can associate the following four operators.

Let $\Lambda_1(\sigma):H^{1/2}(\partial\Omega)\mapsto H^{-1/2}(\partial\Omega)$ be given by
$$\langle \Lambda_1(\sigma)\psi,\varphi\rangle=
\langle\frac{\partial u}{\partial\nu}|_{\partial\Omega},\varphi\rangle=
\int_{\Omega}\nabla u\cdot\nabla\tilde{\varphi},$$
where $\psi$, $\varphi\in H^{1/2}(\partial\Omega)$,
$u$ solves
\eqref{inscrack}-\eqref{voltbdrycond} and $\tilde{\varphi}$ is any
$H^1(\Omega\backslash\sigma)$ function whose trace on $\partial\Omega$ coincides with
$\varphi$.

Let $\mathcal{N}_1(\sigma):{}_0H^{-1/2}(\partial\Omega)\mapsto {}_0H^{1/2}(\partial\Omega)$
be given by
$$\mathcal{N}_1(\sigma)\eta=u|_{\partial\Omega},$$
where $\eta\in{}_0H^{-1/2}(\partial\Omega)$ and $u$ solves
\eqref{inscrack}-\eqref{currentbdrycond}.

Let $\Lambda_2(\sigma):H^{1/2}(\partial\Omega)\mapsto H^{-1/2}(\partial\Omega)$ be given by
$$\langle \Lambda_2(\sigma)\psi,\varphi\rangle=
\langle\frac{\partial u}{\partial\nu}|_{\partial\Omega},\varphi\rangle=
\int_{\Omega}\nabla u\cdot\nabla\tilde{\varphi},$$
where $\psi$, $\varphi\in H^{1/2}(\partial\Omega)$, $u$ solves
\eqref{condcrack}-\eqref{voltbdrycond} and $\tilde{\varphi}$ is any
$H^1_{const}(\Omega,\sigma)$ function whose trace on $\partial\Omega$ coincides
with $\varphi$.

Let $\mathcal{N}_2(\sigma):{}_0H^{-1/2}(\partial\Omega)\mapsto {}_0H^{1/2}(\partial\Omega)$
be given by
$$\mathcal{N}_2(\sigma)\eta=u|_{\partial\Omega},$$
where $\eta\in{}_0H^{-1/2}(\partial\Omega)$ and $u$ solves
\eqref{condcrack}-\eqref{currentbdrycond}.

Let us remark that for any $\sigma\in X$, each $\Lambda_i(\sigma)$ and $\mathcal{N}_i(\sigma)$,
$i=1,2$, is a bounded linear operator between a Hilbert space and its dual, it is self-adjoint
and its norm is always assumed to be the canonical one of bounded operators between these two Hilbert spaces.
Keeping in mind these notations and this remark, we are able to state our instability result.

\begin{prop}\label{invcrackprop}
Let us fix integers $N\geq 2$ and
$m\geq 1$ and a positive constant $\beta$.
Let $(X,d)$ be the metric space where
$X=X_{m\beta(1/4)}(B'_{N-1}(0,1/2))$ and $d$ is the Hausdorff distance.
Let us fix $T\in\{\Lambda_1,\mathcal{N}_1,\Lambda_2,\mathcal{N}_2\}$.
Then there exists a positive $\varepsilon_1$, depending on $N$, $m$ and $\beta$ only,
so that for any $\varepsilon$, $0<\varepsilon<\varepsilon_1$,
there exist two cracks $\sigma_1$,
$\sigma_2$ belonging to $X$
satisfying
\begin{equation}
\begin{array}{l}
d(\sigma_j,\overline{B'_{N-1}(0,1/2)})\leq\varepsilon,\text{ for any }j=1,2;\quad
d(\sigma_1,\sigma_2)\geq\varepsilon;\\
\|T(\sigma_1)-T(\sigma_2)\|\leq2\exp(-\varepsilon^{-(N-1)/(2mN)}).
\end{array}
\end{equation}
\end{prop}

For the proof of this proposition we refer to Subsection~\ref{cracksubproof}.

\subsection{Inverse cavity problem}\label{cavssub}

The inverse
cavity problem can be treated 
if we substitute, in the previous subsection, the
set of cracks inside $\Omega$,
$X=X_{m\beta(1/4)}(B'_{N-1}(0,1/2))$, with the set of \emph{cavities} inside
$\Omega$ given
by
$X=Y_{m\beta(1/4)}(S^{N-1}(0,1/2))$. With almost no modification in the proof,
a result completely analogous to the one
described in Proposition~\ref{invcrackprop} can be obtained. 
So also the inverse cavity problem shows an
exponential instability character.

We recall that, concerning the inverse cavity problem, stability estimates of
logarithmic type have been obtained for the two dimensional case in \cite{Ale e Ron99}
and for the higher dimensional case in \cite{Ale e Bere e Ross e Ves}. For the planar case,
an explicit example developed in \cite{Ale e Ron99} shows the exponential instability
character of the inverse cavity problem and, consequently,
that the stability estimates therein
contained are essentially optimal. Our results here confirm this fact and
extend it to the higher dimensional case, thus providing the essential optimality
of the estimates proved in \cite{Ale e Bere e Ross e Ves}.

\subsection{Inverse surface crack problem}\label{surfcracksub}

Let $\Omega=B_N(0,1)\backslash\{x\in B_N(0,1):\ x_{N-1}\geq 0\text{ and }x_N= 0\}$,
$N\geq 2$,
and let
$\Gamma_I=\{x\in \overline{B_N(0,1)}:\ x_{N-1}\geq 0\text{ and }x_N= 0\}$.
Inside $\Omega$, we consider the geodesic distance between two points as the infimum
of the lengths of smooth paths contained in $\Omega$ connecting the two points.
If we consider the boundary of $\Omega$ with respect to this distance, we notice that
this boundary contains two overlapping copies of $\Gamma_I$, one obtained by approaching
$\Gamma_I$ with points $x$ in $\Omega$ such that $x_N>0$ and the other
obtained by approaching it
with points $x\in\Omega$ so that $x_N<0$.
The set $\Gamma_A$ is obtained from $\Gamma_I$ by taking  the closure, in the topology
of $\partial\Omega$ induced by the geodesic
distance defined above, of $\partial\Omega\backslash \Gamma_I$. We remark that
$\Gamma_A$ coincides, from a set point of view, with $S^{N-1}$, but
each point belonging to the intersection of $\Gamma_I$ and $S^{N-1}$ should be
counted with multiplicity two, as for points of $\Gamma_I$.

With $H^{1/2}(\Gamma_A)$ we denote the space of traces of $H^1(\Omega)$ functions
on $\Gamma_A$ and with $H^{-1/2}(\Gamma_A)$ we shall denote its dual.
On these two spaces, we consider scalar products which are defined exactly
as we have done before for regular domains,
in \eqref{scpr1/2} and \eqref{scpr-1/2}, respectively. We notice that
$H^{1/2}(\Gamma_A)\subset L^2(S^{N-1})\subset
H^{-1/2}(\Gamma_A)$. Finally, we notice that the spaces ${}_0H^{1/2}(\Gamma_A)$ and
${}_0H^{-1/2}(\Gamma_A)$ are the orthogonal subspaces, respectively in
$H^{1/2}(\Gamma_A)$ and
$H^{-1/2}(\Gamma_A)$, to the constant function $1$ and are dual one to each other.

We observe that the spaces $H^1_0(\Omega,\Gamma_I)$ and $H^1_{const}(\Omega,\Gamma_I)$
are given by the spaces of 
$H^1(B_N(0,1))$ functions which are, respectively, identically zero or
constant in a weak sense on $\Gamma_I$. The spaces of traces on $\Gamma_A$ of
functions belonging to $H^1_0(\Omega,\Gamma_I)$ and $H^1_{const}(\Omega,\Gamma_I)$,
respectively, are again denoted by
$H^{1/2}_0(\Gamma_A,\Omega)$ and $H^{1/2}_{const}(\Gamma_A,\Omega)$. On these two last spaces,
a scalar product is defined in the same fashion as we have done for regular
domains in \eqref{scpr1/2const}.

If $N=2$, let $\sigma_0=\{x\in \overline{B_2(0,1)}:\ x_{1}\geq -1/2\text{ and }x_2= 0\}$.
If $N\geq 3$,
let $f\in C^{\infty}_0(B'_{N-2}(0,1/4))$ so that $-1/4\leq f\leq 0$.
Let $\sigma_0=\Gamma_I\cup\{y\in B'_{N-1}(0,1):\
f((y_1,\ldots,y_{N-2},0))\leq y_{N-1}\leq 0,\  (y_1,\ldots,y_{N-2},0)\in
B'_{N-2}(0,1/4)\}$. By definition
if $N=2$, and by a suitable choice of $f$ if $N\geq 3$, 
we can always assume that $B'_{N-1}(\tilde{x}_0,1/16)$ is contained in $\sigma_0$, where
$\tilde{x}_0=(0,\ldots,0,-1/8 ,0)$.

Then we fix a positive integer $m$ and a positive constant $\beta$ and we define
$X$ as the set
\begin{equation}\label{defXsurcrack}
X=\{\sigma=
(\sigma_0\backslash B'_{N-1}(\tilde{x}_0,1/16))\cup\sigma':\ 
\sigma'\in X_{m\beta(1/4)}(B'_{N-1}(\tilde{x}_0,1/16))\}.
\end{equation}

We remark that each $\sigma\in X$ is a connected closed set inside
$\overline{B_N(0,1)}$ so that
$\Gamma_I\subset\sigma$ and $\sigma\backslash\Gamma_I\subset B_N(0,4/5)$.

If we assume that $B_N(0,1)$ is occupied by a homogeneous conductor, we can think
any $\sigma\in X$ as a \emph{surface crack} inside $B_N(0,1)$.
We can distinguish between two different kinds of surface cracks, namely
\emph{insulating} and \emph{conducting}.

Let us assume that $\sigma\in X$ is an insulating surface crack and that we prescribe
on $\Gamma_A$ the current density
to be equal to $\eta\in{}_0H^{-1/2}(\Gamma_A)$. Then
the electrostatic potential $u$ in $B_N(0,1)$ satisfies
\begin{equation}\label{inssurfcrack}
\left\{\begin{array}{ll}
\Delta u=0 &\text{in }B_N(0,1)\backslash\sigma,\\
\frac{\partial u}{\partial\nu}=0 &\text{on either sides of }\sigma,\\
\frac{\partial u}{\partial\nu}=\eta &\text{on }\Gamma_A,\\
\int_{\Gamma_A}u=0,
\end{array}
\right.
\end{equation}
where we have also added a normalization condition. We have that $u$ is a weak solution
to \eqref{inssurfcrack} if and only if
$u\in H^1(B_N(0,1)\backslash\sigma)$, $\int_{\Gamma_A}u=0$, and
$$\int_{B_N(0,1)\backslash\sigma}\nabla u\cdot\nabla w=\eta(w|_{\Gamma_A}),\quad
\text{for any }w\in H^1(B_N(0,1)\backslash\sigma).$$
Clearly such a function $u$ exists and is unique. We have that $u|_{\Gamma_A}$
belongs to ${}_0H^{1/2}(\Gamma_A)$ and that the operator
$\mathcal{N}_3(\sigma):{}_0H^{-1/2}(\Gamma_A)\mapsto{}_0H^{1/2}(\Gamma_A)$ so that,
for any $\eta\in{}_0H^{-1/2}(\Gamma_A)$,
$\mathcal{N}_3(\sigma)\eta=u|_{\Gamma_A}$, $u$ solution to
\eqref{inssurfcrack}, is linear, bounded and self-adjoint.

When, otherwise, $\sigma\in X$ is a conducting surface crack in $B_N(0,1)$
and we prescribe the voltage on $\Gamma_A$
to be $\psi\in H^{1/2}_{const}(\Gamma_A,\Omega)$, then the potential
$u$ in $B_N(0,1)$ solves
\begin{equation}\label{condsurfcrack}
\left\{\begin{array}{ll}
\Delta u=0 &\text{in }B_N(0,1)\backslash\sigma,\\
u=c(\psi) &\text{on }\sigma,\\
u=\psi &\text{on }\Gamma_A,
\end{array}
\right.
\end{equation}
where $c(\psi)$ is a constant so that
$\hat{\psi}=\psi-c(\psi)\in H^{1/2}_0(\Gamma_A,\Omega)$.
Let $\tilde{\psi}$ be any $H^1_0(B_N(0,1)\backslash\sigma,\sigma)$ function so that
$\tilde{\psi}|_{\Gamma_A}=\hat{\psi}$. Then $u$ solves in a weak sense
\eqref{condsurfcrack}
if and only if $u-c(\psi)-\tilde{\psi}\in H^1_0(B_N(0,1)\backslash\sigma)$ and
$$\int_{B_N(0,1)\backslash\sigma}\nabla u\cdot\nabla w=0,\quad
\text{for any }w\in H^1_0(B_N(0,1)\backslash\sigma).$$
By standard elliptic equations methods we infer that $u$, solution to
\eqref{condsurfcrack}, exists and it is unique. To such a solution we can associate
$\frac{\partial u}{\partial \nu}|_{\Gamma_A}
\in(H^{1/2}_{const}(\Gamma_A,\Omega))'$ as follows. For any
$\varphi\in H^{1/2}_{const}(\Gamma_A,\Omega)$, let the constant $c(\varphi)$ be so that
$\hat{\varphi}=\varphi-c(\varphi)\in H^{1/2}_0(\Gamma_A,\Omega)$, and let
$\tilde{\varphi}$ be any $H^1_0(B_N(0,1)\backslash\sigma,\sigma)$ function so that
$\tilde{\varphi}|_{\Gamma_A}=\hat{\varphi}$.
Then,
$$\langle \frac{\partial u}{\partial \nu}|_{\Gamma_A},\varphi\rangle=
\int_{B_N(0,1)\backslash\sigma}\nabla u\cdot\nabla\tilde{\varphi}.$$
The
operator $\Lambda_3(\sigma):H^{1/2}_{const}(\Gamma_A,\Omega)\mapsto (H^{1/2}_{const}(\Gamma_A,\Omega))'$ so that,
for any $\psi\in H^{1/2}_{const}(\Gamma_A,\Omega)$,
$\Lambda_3(\sigma)\psi=\frac{\partial u}{\partial \nu}|_{\Gamma_A}$, $u$ solution to
\eqref{condsurfcrack}, is linear, bounded and self-adjoint.

The inverse surface crack problem consists of the determination of an
unknown surface crack
from suitable information on the operator $\mathcal{N}_3$ or $\Lambda_3$, respectively. The operators
$\mathcal{N}_3$ and $\Lambda_3$ correspond to electrostatic boundary measurements. Many papers have
treated this problem when a finite number of measurements is performed, that is
when either $\mathcal{N}_3(\eta)$ is measured for a finite number of different $\eta$ or
$\Lambda_3(\psi)$ is measured for a finite number of different $\psi$.
We refer to \cite{Ron99:2} and its bibliography for a detailed description
of the problem in the planar case. In \cite{Ron99:2} uniqueness and stability
estimates of logarithmic type are established. In \cite{Ale e DiC} the uniqueness
issue for the three-dimensional case is treated. In the next proposition
we show that also this inverse problem
is exponentially unstable, thus proving the essential optimality of
the stability estimates of \cite{Ron99:2}.

\begin{prop}\label{surfcrackprop}
We fix integers $N\geq 2$ and
$m\geq 1$ and a positive constant $\beta$.
Let $X$ be the set of closed sets described in \eqref{defXsurcrack}
and let
$(X,d)$ be a metric space with the Hausdorff distance. Let $\sigma_0$ be defined as
before.
Let us fix $T\in\{\mathcal{N}_3,\Lambda_3\}$.
Then we can find $\varepsilon_1>0$, that depends on $N$, $m$ and $\beta$ only, so
that for any $\varepsilon$, $0<\varepsilon<\varepsilon_1$,
there exist two surface cracks $\sigma_1$,
$\sigma_2$ belonging to $X$ so that
\begin{equation}
\begin{array}{l}
d(\sigma_j,\sigma_0)\leq\varepsilon,\text{ for any }j=1,2;\quad
d(\sigma_1,\sigma_2)\geq\varepsilon;\\
\|T(\sigma_1)-T(\sigma_2)\|\leq2\exp(-\varepsilon^{-(N-1)/(2mN)}).
\end{array}
\end{equation}
\end{prop}

We prove this result in Subsection~\ref{otherdefsubproof}.

\subsection{Inverse boundary material loss problem}\label{corrsub}

Let $\Omega=\{x\in B_N(0,1):\ x_N> 0\}$, and let
$\Gamma_A=\{x\in \partial B_N(0,1):\ x_N\geq 0\}=S^{N-1}_+$ and
$\Gamma_I=\{x\in \overline{B_N(0,1)}:\ x_N= 0\}=\overline{B'_{N-1}(0,1)}$.
Fixed a positive integer $m$ and a positive constant $\beta$,
let
\begin{equation}\label{defXbml}
X=\{\sigma=\Gamma_I\cup\sigma':\ \sigma'\in Y_{m\beta(1/4)}(B'_{N-1}(0,1/2))\}.
\end{equation}
Then every
$\sigma\in X$ is a closed subset contained in $\overline{\Omega}$ so that
$(\sigma\backslash\Gamma_I)\subset B_N(0,4/5)$.

We assume that $\Omega$ is the region occupied by a homogeneous conductor and
$\sigma\in X$ is a \emph{boundary material loss},
which might be due to a corrosion phenomenon,
for instance. We assume that $\Gamma_A$ is an accessible part of the boundary of
our conductor, whereas
$\Gamma_{\sigma}=\partial(\Omega\backslash\sigma)\backslash\Gamma_A$, that is
the other part of
the boundary where the material loss occurs, is not. 
Also in this case we distinguish two kinds of boundary material losses,
insulating and conducting. In the first case, no current passes through
$\Gamma_{\sigma}$,
the part of
boundary of $\Omega\backslash\sigma$ which is contained in $\sigma$.
In the second case, the voltage is constant on $\sigma$.
More precisely, we have that
if $\sigma$ is an insulating boundary material loss and if we prescribe the current
density on $\Gamma_A$ to be equal to $\eta\in {}_0H^{-1/2}(\Gamma_A)$, then the
electrostatic potential $u$ inside $\Omega\backslash\sigma$ is
the unique solution to
\begin{equation}\label{insbml}
\left\{
\begin{array}{ll}
\Delta u=0 &\text{in }\Omega\backslash\sigma,\\
\frac{\partial u}{\partial\nu}=0&\text{on }\Gamma_{\sigma},\\
\frac{\partial u}{\partial\nu}=\eta&\text{on }\Gamma_A,\\
\int_{\Gamma_A}u=0,
\end{array}
\right.
\end{equation}
where the last line is a normalization condition. Otherwise, if $\sigma$ is conducting,
then the electrostatic potential $u$ in $\Omega$ is given by
\begin{equation}\label{condbml}
\left\{
\begin{array}{ll}
\Delta u=0 &\text{in }\Omega\backslash\sigma,\\
u=c(\psi)&\text{on }\sigma,\\
u=\psi&\text{on }\Gamma_A,
\end{array}
\right.
\end{equation}
where $\psi\in H^{1/2}_{const}(\Gamma_A,\Omega)$ is the prescribed voltage
on $\Gamma_A$ and $c(\psi)$ is a constant so that $\psi-c(\psi)$ belongs to
$H^{1/2}_0(\Gamma_A,\Omega)$.

In the insulating case,
for every $\sigma\in X$, we define $\mathcal{N}_4(\sigma):{}_0H^{-1/2}(\Gamma_A)\mapsto
{}_0H^{1/2}(\Gamma_A)$ so that for any $\eta\in{}_0H^{-1/2}(\Gamma_A)$, then
$\mathcal{N}_4(\sigma)\eta=u|_{\Gamma_A}$, $u$ being the unique solution to \eqref{insbml}.
We have that $\mathcal{N}_4(\sigma)$ is a linear, bounded and
self-adjoint operator.

In the conducting case, if $\sigma\in X$, let us define
$\Lambda_4(\sigma): H^{1/2}_{const}(\Gamma_A,\Omega)\mapsto( H^{1/2}_{const}(\Gamma_A,\Omega))'$ as follows.
For any $\psi$, $\varphi\in H^{1/2}_{const}(\Gamma_A,\Omega)$,
$$\langle \Lambda_4(\sigma)\psi,\varphi\rangle=
\langle \frac{\partial u}{\partial\nu}|_{\Gamma_A},\varphi\rangle=
\int_{\Omega\backslash\sigma}\nabla u\cdot\nabla\tilde{\varphi},$$
where $u$ solves \eqref{condbml} and $\tilde{\varphi}$ is any
$H^1_{const}(\Omega,\sigma)$ so that $\tilde{\varphi}|_{\Gamma_A}=\varphi$.
Also $\Lambda_4(\sigma)$ is a linear, bounded and self-adjoint operator, for any $\sigma\in X$.

The inverse problem consists of the determination of the shape and the location
of an unknown boundary material loss $\sigma$ from electrostatic measurements
performed on the accessible part of the boundary, that is $\Gamma_A$.
The case of a single electrostatic measurement is particularly interesting
and uniqueness and stability estimates have been obtained for this kind of problem,
see \cite{Ron99:2}, and its references, for the two-dimensional case and
\cite{Ale e Bere e Ross e Ves} and also \cite{Che e Hon e Yam}
for the higher dimensional case. The stability estimates obtained in
\cite{Ron99:2} and \cite{Ale e Bere e Ross e Ves} are of
logarithmic type and they are essentially optimal. In two dimensions, this has been shown
through an explicit example provided in \cite{Ale97}. In the next proposition,
proven in Subsection~\ref{otherdefsubproof}, we confirm
that logarithmic stability is essentially optimal in any dimension, no matter how many
and which measurements we perform.

\begin{prop}\label{bmlprop}
Let $N\geq 2$ and
$m\geq 1$ be integers and $\beta$ be a positive constant.
Let $X$ be defined as in \eqref{defXbml}, endowed with the Hausdorff distance $d$.
Let us fix $T\in\{\mathcal{N}_4,\Lambda_4\}$.
Then there exists a constant $\varepsilon_1>0$, that depends on $N$, $m$ and $\beta$ only,
so that for any $\varepsilon$, $0<\varepsilon<\varepsilon_1$,
there exist two boundary material losses $\sigma_1$,
$\sigma_2$ belonging to $X$ so that
\begin{equation}
\begin{array}{l}
d(\sigma_j,\overline{B'_{N-1}(0,1)})\leq\varepsilon,\text{ for any }j=1,2;\quad
d(\sigma_1,\sigma_2)\geq\varepsilon;\\
\|T(\sigma_1)-T(\sigma_2)\|\leq2\exp(-\varepsilon^{-(N-1)/(2mN)}).
\end{array}
\end{equation}
\end{prop}

\subsection{Inverse scattering problem}\label{scatsub}

We turn our attention to inverse scattering problems,
in particular to the determination of \emph{obstacles} inside a medium
by acoustic far-field measurements. For a detailed description
of this kind of inverse problems we refer to \cite{Col e Kre98}.

Let us fix integers $N\in\{2,3\}$ and $m\geq 1$ and two positive constants
$\beta$ and $a$.

Let $X=Y_{m\beta(1/2)}(S^{N-1}(0,1))$. We assume that $D\in X$ is an obstacle in
an otherwise homogeneous medium.

The incident field is determined by a time-harmonic acoustic plane wave and
is given by $u^i(x;\omega;a)=\rme^{\rmi\sqrt{a}x\cdot\omega}$, where $x\in\mathbb{R}^N$,
$\sqrt{a}$ is the wave number, and $\omega\in S^{N-1}$ is the direction of
propagation. The direct scattering problem consists of finding
the total field $u(x;\omega;a)$, $x\in\mathbb{R}^N\backslash D$, which is the sum
of the incident field $u^i(x;\omega;a)$ and of the scattered field
$u^s(x;\omega;a)$, which is due to the presence of the obstacle $D$. The total field $u$
satisfies
\begin{equation}\label{Helmeq}
\left\{\begin{array}{ll}
\Delta u+au=0&\text{in }\mathbb{R}^N\backslash D,\\
u(x;\omega;a)=\rme^{\rmi\sqrt{a}x\cdot\omega}+u^s(x;\omega;a)
&\text{for any }x\in\mathbb{R}^N\backslash D,
\end{array}\right.
\end{equation}
with a boundary condition on $\partial D$ which depends on the nature of the obstacle.
Namely, if the obstacle is \emph{sound-soft}, then
\begin{equation}\label{softbc}
u=0\quad\text{on }\partial D,
\end{equation}
if the obstacle is \emph{sound-hard}, then
\begin{equation}\label{hardbc}
\frac{\partial u}{\partial\nu}=0\quad\text{on }\partial D.
\end{equation}
Furthermore, the scattered field satisfies the so-called Sommerfeld radiation condition
\begin{equation}\label{Somm}
\lim_{r\to\infty}r^{(N-1)/2}\left(\displaystyle{\frac{\partial u^s}{\partial r}}-
\rmi\sqrt{a}u^s\right)=0,
\end{equation}
where $r=\|x\|$ and the limit holds
uniformly for all directions $\hat{x}=x/\|x\|\in S^{N-1}$.
The Sommerfeld radiation
condition characterizes outgoing waves and implies that the
asymptotic behaviour of the scattered field is given by
\begin{equation}\label{asympt}
u^s(x;\omega;a)=\frac{\rme^{\rmi\sqrt{a}\|x\|}}{\|x\|^{(N-1)/2}}\left\{
u^s_{\infty}(\hat{x};\omega;a)+O\left(\frac{1}{\|x\|}\right)
\right\},
\end{equation}
as $\|x\|$ goes to $\infty$, uniformly in all directions $\hat{x}=x/\|x\|\in S^{N-1}$.
The function $u^s_{\infty}$ is called the \emph{far-field pattern} related
to the solution to \eqref{Helmeq}-\eqref{softbc}-\eqref{Somm}, or
to \eqref{Helmeq}-\eqref{hardbc}-\eqref{Somm} respectively, for the direction
of propagation $\omega$ and the wave number $\sqrt{a}$. 

Therefore, if $D\in X$ is sound-soft, then we denote
with $\mathcal{A}_s(D):S^{N-1}\times S^{N-1}\times (0,\infty)\mapsto\mathbb{C}$ its
\emph{far-field map}, that is, for any $\hat{x}$,
$\omega\in S^{N-1}$ and any $a>0$,
\begin{equation}\label{softffpdef}
\mathcal{A}_s(D)(\hat{x},\omega,a)=u^s_{\infty}(\hat{x};\omega;a),
\end{equation}
where $u^s_{\infty}$ is the far-field pattern related
to the solution to \eqref{Helmeq}-\eqref{softbc}-\eqref{Somm}.

In an analogous way, assuming $D\in X$ to be sound-hard, we can associate to $D$ its
far-field map, given by
$\mathcal{A}_h(D):S^{N-1}\times S^{N-1}\times (0,\infty)\mapsto\mathbb{C}$, so that,
for any $(\hat{x},\omega,a)\in S^{N-1}\times S^{N-1}\times (0,\infty)$,
\begin{equation}\label{hardffpdef}
\mathcal{A}_h(D)(\hat{x},\omega,a)=u^s_{\infty}(\hat{x};\omega;a),
\end{equation}
where $u^s_{\infty}$ is now the far-field pattern related
to the solution to \eqref{Helmeq}-\eqref{hardbc}-\eqref{Somm}.

The inverse scattering problem consists of the determination of an unknown
obstacle $D$ from some information about its far-field map. More precisely,
we assume that we a priori know whether the unknown obstacle is sound-soft or sound-hard.
Then, from a, usually partial, knowledge of $\mathcal{A}_s(D)$ or $\mathcal{A}_h(D)$, respectively,
we try to determine the obstacle $D$. We remark that information about the far-field map
can be collected by performing suitable far-field acoustic measurements.

In the next proposition we show that also this inverse problem is
exponentially unstable. We recall that
stability estimates for the determination of a sound-soft obstacle
have been obtained by V.~Isakov, see \cite{Isak92,Isak}.

\begin{prop}\label{scatprop}
Let $N\in\{2,3\}$. Let $m$ and $j$ be positive integers
and $\beta$ be a positive constant. Let $\underline{a}$ and $\overline{a}$
be two
constants so that $0<\underline{a}\leq\overline{a}$. We denote
$$I_N=\left\{\begin{array}{ll}
[\underline{a},\overline{a}] & \text{if }N=2;\\
(0,\overline{a}] & \text{if }N=3;
\end{array}\right.$$
and we fix any $j$ real numbers
$a_1,\ldots,a_j$ so that $a_i\in I_N$
for any $i$, $i=1,\ldots,j$.
Let $(X,d)$ be the metric space with
$X=Y_{m\beta(1/2)}(S^{N-1}(0,1))$ and $d$ the Hausdorff distance.
Then we can find a constant $\varepsilon_1>0$, that depends on $N$, $m$, $j$,
$\beta$ and $I_N$ only,
so that for any $\varepsilon$, $0<\varepsilon<\varepsilon_1$,
there exist obstacles $D_1$, $D_2$, $D_3$ and $D_4$ all belonging to $X$
such that
\begin{equation}
\begin{array}{l}
d(D_j,\overline{B_{N}(0,1)})\leq\varepsilon,\quad\text{for any }j=1,\ldots,4;\\
d(D_1,D_2)\geq\varepsilon;\quad
d(D_3,D_4)\geq\varepsilon;
\end{array}
\end{equation}
and, for the sound-soft case,
\begin{equation}
\sup_{a\in \{a_1,\ldots,a_j\}}
\|(\mathcal{A}_s(D_1)-\mathcal{A}_s(D_2))(\cdot,\cdot,a)\|_{L^2(S^{N-1}\times S^{N-1})}
\leq2\exp(-\varepsilon^{-\frac{N-1}{2mN}});
\end{equation}
and, for the sound-hard case,
\begin{equation}
\sup_{a\in \{a_1,\ldots,a_j\}}
\|(\mathcal{A}_h(D_3)-\mathcal{A}_h(D_4))(\cdot,\cdot,a)\|_{L^2(S^{N-1}\times S^{N-1})}
\leq2\exp(-\varepsilon^{-\frac{N-1}{2mN}}).
\end{equation}
\end{prop}

We observe that the result is slightly different depending whether
$N=2$ or $N=3$. If $N=2$, the real numbers $a_i$, $i=1,\ldots,j$,
satisfy $a_i\geq \underline{a}$ and the
constant $\varepsilon_1$ depends on $\underline{a}$ as well. If $N=3$, we do not need
a lower bound for the numbers $a_i$, apart from the fact that they are positive,
and the constant $\varepsilon_1$ does not depend on $\underline{a}$. We shall point out
during the proof, which will be developed in Subsection~\ref{scatsubproof},
where the hypothesis $a_i\geq \underline{a}$ becomes necessary.

\section{The abstract theorem}\label{abstractsec}

Let $(X,d)$ be a metric space and
let $H$ be a separable Hilbert space with scalar product $(\cdot,\cdot)$.
As usual we denote with $H'$ its dual and for any $v'\in H'$ and any $v\in H$
we denote by $\langle v',v \rangle$ the duality pairing between $v'$ and $v$.
With $\mathcal{L}(H,H')$ we denote the space of bounded linear operators between
$H$ and $H'$ with the usual operators norm.
We shall also fix
$\gamma:H\backslash\{0\}\mapsto [0,+\infty]$ such that
\begin{equation}
\gamma(\lambda v)=\gamma(v)\quad\text{for any }v\in
H\backslash\{0\}\text{ and any }\lambda\in \mathbb{R}\backslash\{0\}.
\end{equation}
Let us remark that the function $\gamma$ may attain both the values $0$ and
$+\infty$ and can be thought of as a suitable kind of
Rayleigh quotient.

Let $F$ be a function from $X$ to $\mathcal{L}(H,H')$, that is, for any $x\in X$,
$F(x)$ will denote a linear and bounded operator between $H$ and $H'$. We also fix a
reference operator $F_0\in\mathcal{L}(H,H')$ and a reference point $x_0$ in $X$.
We wish to point out that no connection is required between $x_0$ and $F_0$, in
particular $F_0$ might be different from $F(x_0)$.
For any $\varepsilon>0$, let $X_{\varepsilon}=\{x\in X:\ d(x,x_0)\leq\varepsilon\}$. 

Recalling the notations introduced in Definition~\ref{def-net-discrete},
we can formulate the following exponential instability result related to
the map $F$.

\begin{teo}\label{abstractthm}
Let us assume that the following conditions are satisfied.
\begin{enumerate}[i\textnormal{)}]
\item\label{Xassumpt}
There exist positive constants $\varepsilon_0$, $C_1$ and $\alpha_1$ such that
for any $\varepsilon$, $0<\varepsilon<\varepsilon_0$, we can find an
$\varepsilon$-discrete set $Z_{\varepsilon}$ contained in $X_{\varepsilon}$ with at
least $\exp(C_1\varepsilon^{-\alpha_1})$ elements.

\item\label{HFassumpt}
There exist three positive constants $p$,
$C_2$ and $\alpha_2$
and
an orthonormal basis in $H$, $\{v_k\}_{k=1}^{+\infty}$, such that the following
conditions hold.

For any $k\in\mathbb{N}$, we have that $\gamma(v_k)<\infty$, and
for any $n\in\mathbb{N}$,
\begin{equation}\label{polygrowth}
\#\{k\in\mathbb{N}:\ \gamma(v_k)\leq n\}\leq C_2(1+n)^{p}
\end{equation}
where $\#$ denotes the number of elements.

For any $x\in X$ and any $(k,l)\in\mathbb{N}\times\mathbb{N}$
we have
\begin{equation}\label{mainestimate}
|\langle (F(x)-F_0)v_k,v_l\rangle|\leq
C_2\exp(-\alpha_2\max\{\gamma(v_k),\gamma(v_l)\}).
\end{equation}
\end{enumerate}

Then there exists a positive constant $\varepsilon_1$,
depending on $\varepsilon_0$, $C_1$, $C_2$, $\alpha_1$, $\alpha_2$ and $p$ only,
so that for every $\varepsilon$, $0<\varepsilon<\varepsilon_1$, we can find
$x_1$ and $x_2$ satisfying
\begin{equation}
\begin{array}{l}
x_1,x_2\in X_{\varepsilon};\quad
d(x_1,x_2)\geq\varepsilon;\\
\|F(x_1)-F(x_2)\|_{\mathcal{L}(H,H')}\leq 2\exp(-\varepsilon^{-\alpha_1/2(p+1)}).
\end{array}
\end{equation}
\end{teo}

A crucial point in the proof of Theorem~\ref{abstractthm} is constituted by the
following lemma in which we construct $\delta$-nets in the image through $F$ of
$X$ with a control in terms of $\delta$ of the number of their elements.

\begin{lem}\label{netlemma}
Under assumption \ref{HFassumpt}\textnormal{)}
of Theorem~\textnormal{\ref{abstractthm}},
there exists a positive constant $C_3$,
depending on $p$, $C_2$ and $\alpha_2$ only,
such that for every $\delta$, $0<\delta<1/\rme$, we can find a $\delta$-net
$Y_{\delta}$ for $F(X)$ with at most
$\exp(C_3(-\log\delta)^{2p+1})$ elements.
\end{lem}

\proof{.} Let $G$ be an element of $\mathcal{L}(H,H')$ and let 
$\{v_k\}_{k=1}^{+\infty}$ be the orthonormal basis in $H$ defined in
assumption \textit{\ref{HFassumpt}}\textnormal{)} of
Theorem~\textnormal{\ref{abstractthm}}.
For any pair $(k,l)\in\mathbb{N}\times\mathbb{N}$ we associate to $G$ the
real number $a_{k,l}=\langle Gv_k,v_l \rangle$. Let
$\|G\|_{Y}=\sup_{k,l}|a_{k,l}|(2+\max\{\gamma(v_k),\gamma(v_l)\})^{p+1}$
and let $Y$ be the normed space
$$Y=\{G\in\mathcal{L}(H,H'):\  \|G\|_{Y}<\infty\}$$
with norm $\|\cdot\|_{Y}$.
First, we notice that, for any $x\in X$, $F(x)-F_0$ is contained in $Y$.
This is an immediate consequence of \eqref{mainestimate}; in fact
$|\langle (F(x)-F_0)v_k,v_l\rangle|\leq
C_2\exp(-\alpha_2\max\{\gamma(v_k),\gamma(v_l)\})$ and hence
$$\|F(x)-F_0\|_Y\leq \sup_n C_2\exp(-\alpha_2(n-1))(2+n)^{p+1}<\infty.$$
Second, if we set $C_4=C_2\left(\sum_n (1+n)^{-2}\right)^{1/2}$, for any $G\in Y$ we have
\begin{equation}\label{compar}
\|G\|_{\mathcal{L}(H,H')}\leq C_4\|G\|_{Y}.
\end{equation}
This follows from the following computation
$$\begin{array}{rl}
\|G\|_{\mathcal{L}(H,H')}&\leq\left(\sum_{k,l}|a_{k,l}|^2\right)^{1/2}\leq\\
&\leq\left(\sum_{k,l}|a_{k,l}|^2\frac{(2+
\max\{\gamma(v_k),\gamma(v_l)\})^{2p+2}}{(2+
\max\{\gamma(v_k),\gamma(v_l)\})^{2p+2}}\right)^{1/2}\leq\\
&\leq\left(\sum_{k,l}\frac{1}{(2+
\max\{\gamma(v_k),\gamma(v_l)\})^{2p+2}}\right)^{1/2}\|G\|_{Y}.
\end{array}$$
Let us show that $C_4$ is a bound for 
$\left(\sum_{k,l}\frac{1}{(2+
\max\{\gamma(v_k),\gamma(v_l)\})^{2p+2}}\right)^{1/2}$.
For any positive integer $n$, the number of pairs $(k,l)$ so that
$n-1\leq\max\{\gamma(v_k),\gamma(v_l)\}< n$ is bounded by $C_2^2(1+n)^{2p}$, therefore
$$\sum_{k,l}\frac{1}{(2+
\max\{\gamma(v_k),\gamma(v_l)\})^{2p+2}}\leq
\sum_n C_2^2(1+n)^{2p}\frac{1}{(2+(n-1))^{2p+2}}=C_4^2.$$

Let us now fix $\delta$, $0<\delta<1/\rme$. Let $\tilde{n}$ be the smallest positive
integer so that for any real number $t\geq\tilde{n}$ it holds
$C_2\exp(-\alpha_2(t-1))(2+t)^{p+1}\leq
\delta/(2C_4)$. There exists a positive constant
$C_5$, depending on $p$, $C_2$ and $\alpha_2$ only, such that
$\tilde{n}\leq C_5\log\delta^{-1}$.

Let $\delta'=\frac{(2+\tilde{n})^{-(p+1)}}{2C_4}\delta$ and let
$\Psi_{\delta}=\delta'\mathbb{Z}\cap[-C_2,C_2]$. We remark that $\Psi_{\delta}$ is a
finite subset of $\mathbb{R}$ and we have that $\#\Psi_{\delta}\leq C_6/\delta'$,
where $C_6$ can be chosen as $2C_2+1$.

Let us define the following subset of $\mathcal{L}(H,H')$. Let
$\tilde{Y}_{\delta}=\{G\in\mathcal{L}(H,H'):\ a_{k,l}\in\Psi_{\delta}\text{ if }
\max\{\gamma(v_k),\gamma(v_l)\}\leq\tilde{n}\text{ and }a_{k,l}=0\text{ otherwise}\}$.

Let us count the number of elements of $\tilde{Y}_{\delta}$. If we set
$$s=\#\{(k,l):\ \max\{\gamma(v_k),\gamma(v_l)\}\leq\tilde{n}\}$$
then we have that
$\#\tilde{Y}_{\delta}=(\#\Psi_{\delta})^{s}$ and hence
$$\begin{array}{rl}
\#\tilde{Y}_{\delta}&\leq(2C_6
C_4(2+\tilde{n})^{p+1}\delta^{-1})^{C_2^2(1+\tilde{n})^{2p}}\leq\\
&\leq(2C_6
C_4(2+C_5\log\delta^{-1})^{p+1}\delta^{-1})^{C_2^2(1+C_5\log\delta^{-1})^{2p}}
\end{array}$$
and then by straightforward computations we can find a positive constant
$C_3$,
depending on $p$, $C_2$ and $\alpha_2$ only, so that
$\#\tilde{Y}_{\delta}\leq\exp(C_3(-\log\delta)^{2p+1})$.

Then we need to show that for every $x\in X$ there exists $G\in\tilde{Y}_{\delta}$
so that $\|F(x)-F_0-G\|_{\mathcal{L}(H,H')}\leq \delta/2$.

We fix $x\in X$ and we set
$b_{k,l}=\langle (F(x)-F_0)v_k,v_l \rangle$; we recall that $F(X)-F_0\subset Y$ and,
by \eqref{compar},
it is enough to determine $G\in\tilde{Y}_{\delta}$ so that
$\|F(x)-F_0-G\|_Y\leq \delta/(2C_4)$.
We observe that, by \eqref{mainestimate}, $b_{k,l}\in[-C_2,C_2]$ for every $(k,l)$.
We construct such a $G$ as follows. We set
$a_{k,l}=\langle Gv_k,v_l \rangle$.
For any $(k,l)$ so that $\max\{\gamma(v_k),\gamma(v_l)\}\leq\tilde{n}$,
we prescribe $a_{k,l}$ to be the element of $\Psi_{\delta}$ that is
closest to $b_{k,l}$. If, otherwise, $(k,l)$ is so that
$\max\{\gamma(v_k),\gamma(v_l)\}>\tilde{n}$, then we set
$a_{k,l}=0$. We have that $G$ belongs to $\tilde{Y}_{\delta}$
by construction and we can evaluate $\|F(x)-F_0-G\|_Y$ as follows. For any $(k,l)$
so that $\max\{\gamma(v_k),\gamma(v_l)\}\leq\tilde{n}$,
$|a_{k,l}-b_{k,l}|\leq\delta'$, that is
$|a_{k,l}-b_{k,l}|(2+\max\{\gamma(v_k),\gamma(v_l)\})^{p+1}\leq
(2+\tilde{n})^{p+1}
\frac{(2+\tilde{n})^{-(p+1)}}{2C_4}\delta\leq\delta/(2C_4)$.
If $(k,l)$ is such that $t=\max\{\gamma(v_k),\gamma(v_l)\}>\tilde{n}$,
then $a_{k,l}=0$ and
$|b_{k,l}|(2+t)^{p+1}\leq
C_2\exp(-\alpha_2t)
(2+t)^{p+1}\leq \delta/(2C_4)$
by the definition of $\tilde{n}$.

Having established this property,
it is easy to find a subset $Y_{\delta}$ of $F(X)$
with
the same number of elements as $\tilde{Y}_{\delta}$ which is a
$\delta$-net for $F(X)$ and hence the proof is concluded.\cvd{}
\proof{ of Theorem~\textnormal{\ref{abstractthm}}.} The proof is obtained
by combining assumption \textit{\ref{Xassumpt}}\textnormal{)} of
Theorem~\ref{abstractthm}
with Lemma~\ref{netlemma} as follows. Let $\varepsilon$ satisfy
$0<\varepsilon<\varepsilon_0$ and $\delta$ satisfy $0<\delta<1/\rme$.
Let $Z_{\varepsilon}\subset X_{\varepsilon}$ be, as in assumption
\textit{\ref{Xassumpt}}\textnormal{)} of Theorem~\ref{abstractthm}, an
$\varepsilon$-discrete set with at
least $\exp(C_1\varepsilon^{-\alpha_1})$ elements. The same procedure employed in
Lemma~\ref{netlemma} allows us to
find $Y_{\delta}\subset F(X_{\varepsilon})$ which is
a $\delta$-net for $F(X_{\varepsilon})$ with at most
$\exp(C_3(-\log\delta)^{2p+1})$ elements.
If $\# Z_{\varepsilon}>\# Y_{\delta}$, then there exist $x_1$ and
$x_2\in X_{\varepsilon}$ so that $d(x_1,x_2)\geq\varepsilon$ and
$\|F(x_1)-F(x_2)\|_{\mathcal{L}(H,H')}\leq 2\delta$.

In order to have that $\# Z_{\varepsilon}>\# Y_{\delta}$, it suffices to have that
$\exp(C_1\varepsilon^{-\alpha_1})>\exp(C_3(-\log\delta)^{2p+1})$. Let us define
$\delta(\varepsilon)=\exp(-\varepsilon^{-\alpha_1/2(p+1)})$.
Then there exists a constant $\varepsilon_1$, $0<\varepsilon_1\leq\varepsilon_0$,
depending on $\varepsilon_0$, $C_1$, $C_2$, $\alpha_1$, $\alpha_2$ and $p$ only,
such that for every $\varepsilon$,
$0<\varepsilon<\varepsilon_1$, we have
$\delta(\varepsilon)<1/\rme$ and
$\exp(C_1\varepsilon^{-\alpha_1})>\exp(C_3(-\log\delta(\varepsilon))^{2p+1})$
and so the result follows.\cvd

\begin{oss}\label{remark1}
We wish to remark that it is easy to see that
the order of instability can be improved up to
$\exp(-\varepsilon^{-\alpha_1/(2p+1+\alpha_3)})$, for any $\alpha_3>0$. However,
in this case,
the constant $\varepsilon_1$ depends on $\alpha_3$, too.
For the sake of simplicity, we have
stated the theorem when $\alpha_3$ is chosen to be equal to $1$.
\end{oss}

\section{Stekloff eigenvalue problems}\label{prelsec}

In this section we collect some results which will be repeatedly used later,
when we shall apply the abstract theorem to find instability examples for
inverse problems. Most of the results described in this section are obtained by standard
methods, thus, for the sake of brevity, we do not enter into any detail and we limit
ourselves to fix the notation and to state the results which will be needed later,
referring to the literature when necessary.

Let $\Omega\subset\mathbb{R}^N$, $N\geq 2$, be a bounded domain
and let $\Gamma_A$
and $\Gamma_I$ be two internally disjoint subsets of $\partial\Omega$,
so that $\Gamma_A\cup\Gamma_I=\partial\Omega$.
About the regularity and the properties of $\Omega$, $\Gamma_A$ and
$\Gamma_I$, we shall consider the same assumptions used at the beginning of
Section~\ref{ressec}, at page~\pageref{sobnot}.

The following eigenvalue problems of Stekloff type will be discussed; first
\begin{equation}\label{PA1}\tag{EP1}
\left\{\begin{array}{ll}
\Delta u=0 & \text{in }\Omega,\\
\frac{\partial u}{\partial \nu}=\lambda u & \text{on }\Gamma_A,\\
\frac{\partial u}{\partial \nu}=0 & \text{on }\Gamma_I,
\end{array}\right.
\end{equation}
and then, assuming $\Gamma_I$ not empty,
\begin{equation}\label{PA2}\tag{EP2}
\left\{\begin{array}{ll}
\Delta v=0 & \text{in }\Omega,\\
\frac{\partial v}{\partial \nu}=\mu v & \text{on }\Gamma_A,\\
v=0 & \text{on }\Gamma_I.
\end{array}\right.
\end{equation}

We state the following propositions concerning the
eigenvalues and eigenfunctions of \eqref{PA1} and \eqref{PA2} respectively.

\begin{prop}\label{PA1prop}
Under the assumptions on $\Omega$, $\Gamma_A$ and $\Gamma_I$ previously made, we have that
the eigenvalues of \eqref{PA1}, counted with their multiplicity,
are given by an increasing sequence
$$0=\lambda_1<\lambda_2\leq\lambda_3\leq\ldots\leq\lambda_k\leq\ldots$$
so that $\lim_{k\to\infty}\lambda_k=\infty$.
For any $n\in\mathbb{N}$, we set $N_1(n)=\#\{k\in\mathbb{N}:\ \lambda_k\leq n\}$. Then
the asymptotic behaviour of the eigenvalues is as follows. There exists a constant $C_1$
depending on $\Omega$, $\Gamma_A$ and $\Gamma_I$ only so that
\begin{equation}\label{eigenasympt1}
N_1(n)\leq C_1n^{N-1},\quad\text{for any }n\in\mathbb{N}.
\end{equation}

Moreover, there exists a corresponding sequence of
eigenfunctions, $\{u_k\}_{k\in\mathbb{N}}$, that is $u_k\in H^1(\Omega)\backslash\{0\}$
and the couple $(\lambda_k,u_k)$ solves
\eqref{PA1} for any $k\in\mathbb{N}$, so that the following three conditions holds
$$\{u_k|_{\Gamma_A}\}_{k\in\mathbb{N}}\text{ is an orthonormal basis of }L^2(\Gamma_A);$$
$$\left\{\frac{u_k}{\sqrt{1+\lambda_k}}|_{\Gamma_A}\right\}_{k\in\mathbb{N}}\text{ is an
orthonormal basis of }H^{1/2}(\Gamma_A);$$
$$\{1|_{\Gamma_A}\}
\cup\{\sqrt{\lambda_k}u_k|_{\Gamma_A}\}_{k\geq 2}\text{ is an
orthonormal basis of }H^{-1/2}(\Gamma_A);$$
where we have considered the spaces $H^{1/2}(\Gamma_A)$ and $H^{-1/2}(\Gamma_A)$ with the
scalar products defined in \eqref{scpr1/2} and \eqref{scpr-1/2} respectively.
We remark that $u_1$ is a constant function not identically
equal to zero.
\end{prop} 
\begin{prop}\label{PA2prop}
Under the assumptions on $\Omega$, $\Gamma_A$ and $\Gamma_I$ previously made,
and assuming that $\Gamma_I$ is not empty, then
the eigenvalues of \eqref{PA2}, counted with their multiplicity,
constitute an increasing sequence
$$0<\mu_1\leq\mu_2\leq\ldots\leq\mu_k\leq\ldots$$
so that $\lim_{k\to\infty}\mu_k=\infty$.
For any $n\in\mathbb{N}$, we set as before $N_2(n)=\#\{k\in\mathbb{N}:\ \mu_k\leq n\}$.
Then the eigenvalues satisfy the following
asymptotic behaviour. There exists a constant $C_2$
depending on $\Omega$, $\Gamma_A$ and $\Gamma_I$ only so that
\begin{equation}\label{eigenasympt2}
N_2(n)\leq C_2n^{N-1},\quad\text{for any }n\in\mathbb{N}.
\end{equation}

Furthermore, we can find a sequence $\{v_k\}_{k\in\mathbb{N}}$ of corresponding
eigenfunctions, that is $v_k\in H^1(\Omega)\backslash\{0\}$
and the couple $(\mu_k,v_k)$ is a solution to
\eqref{PA2} for any $k\in\mathbb{N}$, so that
$$\{v_k|_{\Gamma_A}\}_{k\in\mathbb{N}}\text{ is an orthonormal system of }L^2(\Gamma_A);$$
$$\{1|_{\Gamma_A}\}
\cup\left\{\frac{v_k}{\sqrt{\mu_k}}|_{\Gamma_A}\right\}_{k\in\mathbb{N}}\text{ is an
orthonormal basis of }H^{1/2}_{const}(\Gamma_A,\Omega);$$
where we have considered the space $H^{1/2}_{const}(\Gamma_A,\Omega)$
with the scalar product defined in \eqref{scpr1/2const}.
\end{prop}

Beyond the asymptotic behaviour of the eigenvalues, we are interested
in the asymptotic behaviour of the eigenfunctions, in particular in a kind of exponential
decay, in terms of the eigenvalues, of the eigenfunctions away from $\Gamma_A$.
In the next examples, we present some particular cases in which
such kind of decay holds.

\begin{exam}\label{ballex}
Let $\Omega=B_N(0,1)$ and $\partial\Omega=S^{N-1}$, and let $\Gamma_A=\partial\Omega$
and $\Gamma_I=\emptyset$. In this case the problem \eqref{PA1} is a classical
Stekloff eigenvalue problem and it is well-known that the orthonormal basis of
$L^2(S^{N-1})$ constituted by the traces of eigenfunctions, as described in
Proposition~\ref{PA1prop}, coincides with
\begin{equation}\label{spherarm}
\{f_{jp}:\ j\geq 0\text{ and }1\leq p \leq p_j\}
\end{equation}
where each $f_{jp}$ is a \emph{spherical harmonic} of degree $j$, $j$ being a
nonnegative integer. We have that the function
\begin{equation}\label{polyharm}
u_{jp}(x)=\|x\|^jf_{jp}(x/\|x\|)
\end{equation}
is harmonic in $\mathbb{R}^N$ and solves the eigenvalue problem \eqref{PA1} with
eigenvalue $\lambda=j$. So, the sequence
$\{u_{jp}:\ j\geq 0\text{ and }1\leq p \leq p_j\}$ coincides with the sequence of
eigenfunctions we have described in Proposition~\ref{PA1prop}. The integers $p_j$
are the dimensions of the spaces of spherical harmonics of degree $j$ and we have that,
see for instance \cite[page~4]{Mul},
$$p_j=\left\{\begin{array}{ll}1 & \text{if }j=0,\\
\frac{(2j+N-2)(j+N-3)!}{j!(N-2)!} & \text{if }j\geq 1,\end{array}\right.$$
so that
$$p_j\leq 2(j+1)^{N-2},\quad j\geq 0,$$
and
$$N_1(n)\leq\sum_{j=0}^np_j\leq\sum_{j=0}^n2(j+1)^{N-2}\leq 2(n+1)^{N-1},\quad\text{for
any }n\in\mathbb{N}.$$

Furthermore, for any $r_0$, $0<r_0<1$, there exist two positive constants,
$C(r_0,N)$ and $\alpha(r_0)$, so that for any $u_{jp}$ as in \eqref{polyharm} it holds
\begin{equation}\label{expdecayball}
\|u_{jp}\|_{H^1(B_N(0,r_0))}\leq C(r_0,N)\exp(-\alpha(r_0)j).
\end{equation}
\end{exam}

\begin{exam}\label{halfballex}
Let $\Omega=\{x\in B_N(0,1):\ x_N> 0\}$, and let
$\Gamma_A=\{x\in \partial B_N(0,1):\ x_N\geq 0\}=S^{N-1}_+$ and
$\Gamma_I=\{x\in \overline{B_N(0,1)}:\ x_N= 0\}=\overline{B'_{N-1}(0,1)}$.
First of all, we notice that the hypoteses of Proposition~\ref{PA1prop} and
Proposition~\ref{PA2prop} are satisfied, so the conclusions of
Proposition~\ref{PA1prop} and of Proposition~\ref{PA2prop} hold
for the eigenvalues and eigenfunctions related to
problem \eqref{PA1} and problem \eqref{PA2} with these data, respectively.

The following exponential decay property can be obtained, as well.

We have that if $u\in H^1(\Omega)\backslash\{0\}$ solves
\eqref{PA1} for a constant $\lambda$,
then, by a reflection argument, it follows that there exist $j$, a nonnegative integer,
and $f$, a spherical harmonic function on $S^{N-1}$ of degree $j$, so that
$u(x)=\|x\|^jf(x/\|x\|)$
for any $x\in\Omega$ and $\lambda=j$. Thus, if we assume that
$\|f\|_{L^2(\Gamma_A)}=1$, we can conclude that for any $r_0$, $0<r_0<1$,
\begin{equation}\label{expdecayhalfball1}
\|u\|_{H^1(B_N(0,r_0)\cap\Omega)}\leq C(r_0,N)\exp(-\alpha(r_0)\lambda),
\end{equation}
where the constants $C(r_0,N)$ and $\alpha(r_0)$ coincide with the ones obtained
in Example~\ref{ballex}.

Again by a reflection argument, we have that
if $v\in H^1(\Omega)\backslash\{0\}$ and a constant $\mu$ solve
\eqref{PA2} then there exist $j$, a positive integer,
and $f$, a spherical harmonic function on $S^{N-1}$ of degree $j$, so that
$v(x)=\|x\|^jf(x/\|x\|)$
for any $x\in\Omega$ and $\mu=j$. Thus, if we assume as before that
$\|f\|_{L^2(\Gamma_A)}=1$, we immediately infer that for any $r_0$, $0<r_0<1$,
\begin{equation}\label{expdecayhalfball2}
\|v\|_{H^1(B_N(0,r_0)\cap\Omega)}\leq C(r_0,N)\exp(-\alpha(r_0)\mu),
\end{equation}
with the same constants $C(r_0,N)$ and $\alpha(r_0)$ as before.
\end{exam}

\begin{exam}\label{cutballex}
Let $\Omega=B_N(0,1)\backslash\{x\in B_N(0,1):\ x_{N-1}\geq 0\text{ and }x_N= 0\}$,
and let
$\Gamma_I=\{x\in \overline{B_N(0,1)}:\ x_{N-1}\geq 0\text{ and }x_N= 0\}$,
as in Subsection~\ref{surfcracksub}. We set $\Gamma_A$ as in Subsection~\ref{surfcracksub}, as well.
Also the notations concerning Sobolev spaces on $\Gamma_A$ are the ones
introduced in Subsection~\ref{surfcracksub}.

The eigenvalue problem \eqref{PA1} with these data can be rewritten as
\begin{equation}\label{PA1'}\tag{EP1'}
\left\{\begin{array}{ll}
\Delta u=0 & \text{in }\Omega,\\
\frac{\partial u}{\partial \nu}=\lambda u & \text{on }\Gamma_A,\\
\frac{\partial u}{\partial \nu}=0 & \text{on (either sides of) }\Gamma_I,
\end{array}\right.
\end{equation}
that is, $u\in H^1(\Omega)$ solves \eqref{PA1'} if
$$\int_{\Omega}\nabla u\cdot\nabla w=\int_{\Gamma_A}\lambda uw,\quad\text{for
any }w\in H^1(\Omega).$$

Then, all the conclusions of Proposition~\ref{PA1prop} still hold true for the
eigenvalue problem \eqref{PA1'}, also with the possibility to replace
the space $L^2(\Gamma_A)$ with the space $L^2(S^{N-1})$.

The exponential decay of the eigenfunctions is still valid.
By separation of variables,
we have that if $u\in H^1(\Omega)\backslash\{0\}$ solves \eqref{PA1'} with a constant
$\lambda$, then there exists a function $g\in L^2(S^{N-1})$ so that
$u(x)=\|x\|^{\lambda}g(x/\|x\|)$
for any $x\in\Omega$. Assuming that
$\|g\|_{L^2(S^{N-1})}=1$, we obtain that for any $r_0$, $0<r_0<1$,
\begin{equation}\label{expdecaycutball1}
\|u\|_{H^1(B_N(0,r_0)\cap\Omega)}\leq C_1(r_0,N)\exp(-\alpha(r_0)\lambda),
\end{equation}
where $C_1(r_0,N)$ is a positive constant not depending on $\lambda$ and
$\alpha(r_0)$ coincides with the one defined in Example~\ref{ballex}.

For what concerns the eigenvalue problem \eqref{PA2} with these data, that is,
\begin{equation}\label{PA2'}\tag{EP2'}
\left\{\begin{array}{ll}
\Delta v=0 & \text{in }\Omega,\\
\frac{\partial v}{\partial \nu}=\mu v & \text{on }\Gamma_A,\\
v=0 & \text{on }\Gamma_I,
\end{array}\right.
\end{equation}
we have that $v$ solves \eqref{PA2'}, in a weak sense, for a constant $\mu$, if
$v\in H^1_0(\Omega,\Gamma_I)$ and
$$\int_{\Omega}\nabla v\cdot\nabla w=\int_{\Gamma_A}\mu vw,\quad\text{for
any }w\in H^1_0(\Omega,\Gamma_I).$$

Then, all the results of Proposition~\ref{PA2prop} are still valid for the
eigenvalue problem \eqref{PA2'}, and we can again replace
the space $L^2(\Gamma_A)$ with the space $L^2(S^{N-1})$.

If $v\in H^1_0(\Omega,\Gamma_I)\backslash\{0\}$ and $\mu$ are
a solution to \eqref{PA2'}, then, by separation of variables, we can find
a function $g\in L^2(S^{N-1})$ so that
$v(x)=\|x\|^{\mu}g(x/\|x\|)$
for any $x\in\Omega$. If we further suppose
$\|g\|_{L^2(S^{N-1})}=1$, we have that for any $r_0$, $0<r_0<1$,
\begin{equation}\label{expdecaycutball2}
\|v\|_{H^1(B_N(0,r_0)\cap\Omega)}\leq C_2(r_0,N)\exp(-\alpha(r_0)\mu),
\end{equation}
where $C_2(r_0,N)$ is a positive constant not depending on $\mu$ and
$\alpha(r_0)$ is the same as before.
\end{exam}

\section{Proofs of the main results}\label{examsecproof}

In this section we apply the abstract theorem to the inverse problems
described in Section~\ref{ressec} and we conclude the proofs of our instability results. 

\subsection{Inverse inclusion problem}\label{inclsubproof}

We treat the problem related to the
Dirichlet-to-Neumann map, the one related to the Neumann-to-Dirichlet
map and the one related to the experimental measurements separately.

\subsubsection{Dirichlet-to-Neumann case}

\proof{ of Proposition~\ref{inclprop} (Dirichlet-to-Neumann case).}
We apply Theorem~\ref{abstractthm} to prove Proposition~\ref{inclprop}
for the operator $\Lambda$. Let us show how the abstract theorem is used in this situation.

We fix $x_0\in X$ to be equal to $\overline{B_{N}(0,1/2)}$ and we notice 
that Proposition~\ref{discretesetprop} implies that
$X$ satisfies assumption~\textit{\ref{Xassumpt}}\textnormal{)}
of Theorem~\textnormal{\ref{abstractthm}}, with constants $\varepsilon_0$ and $C_1$
which depend
on $N$, $m$ and $\beta$ only, and with the constant $\alpha_1=(N-1)/m$.
Furthermore, about $X$, we remark the following property. We have
that every $D\in X$ is compactly contained in $B_N(0,4/5)$.

Concerning assumption~\textit{\ref{HFassumpt}}\textnormal{)}
of Theorem~\textnormal{\ref{abstractthm}}, we observe that
for every $D\in X$, the operator $\Lambda(D)$ is a bounded and linear operator between
$H^{1/2}(\partial\Omega)$ and its dual. We fix $H=H^{1/2}(\partial\Omega)$
and $F:X\mapsto\mathcal{L}(H,H')$ as $F(D)=\Lambda(D)$ for any $D\in X$.
With $F_0$ we denote the Dirichlet-to-Neumann map associated to the problem
\eqref{incleq}-\eqref{Dircond} 
when $D=\emptyset$, that is when no inclusion is present, the conductor is
therefore homogeneous and its
conductivity is identically equal to $1$ in $\Omega$.
Concerning the function $\gamma$ our choice is the following.
For any $\psi\in H^{1/2}(\partial\Omega)\backslash\{0\}$, let
\begin{equation}\label{defgammaDir}
\gamma(\psi)=
\frac{\|\psi\|^2_{H^{1/2}(\partial\Omega)}}{\|\psi\|^2_{L^2(\partial\Omega)}}.
\end{equation}

Then it remains to choose an orthonormal basis of $H$, $\{v_k\}_{k\in\mathbb{N}}$,
so that $\gamma(v_k)$ is finite for any $k$ and \eqref{polygrowth}
and \eqref{mainestimate} are satisfied. Recalling Example~\ref{ballex}, in particular \eqref{spherarm},
we consider the set
\begin{equation}\label{basisdefDir}
\left\{\frac{f_{jp}}{\sqrt{1+j}}:\ j\geq 0\text{ and }1\leq p \leq p_j\right\}
\end{equation}
with the natural order. This set, by Proposition~\ref{PA1prop}, is an orthonormal basis
of $H$ and it is the one we choose. We also recall that
$f_{jp}$ is a spherical harmonic of degree $j$ so that
$\|f_{jp}\|_{L^2(\partial\Omega)}=1$, hence $\gamma(f_{jp}/\sqrt{1+j})=1+j$,
for any $j$ and $p$.
Fixed $n\in\mathbb{N}$, $\#\{k\in\mathbb{N}:\ \gamma(v_k)\leq n\}$ is clearly bounded from
above by $2(1+n)^{N-1}$, see Example~\ref{ballex}.

We now verify that \eqref{mainestimate} is also satisfied.
First of all, we remark that for any $D\in X$ we have that the operator
$F(D)-F_0$ is self-adjoint, in the following sense:
$$\langle(F(D)-F_0)\psi,\varphi\rangle=\langle(F(D)-F_0)\varphi,\psi\rangle,$$
for any $\psi$, $\varphi\in H^{1/2}(\partial\Omega)$, where $\langle\cdot,\cdot\rangle$
is as usual the
duality pairing between $H'$ and $H$. In fact, let $u$
be the solution to \eqref{incleq}-\eqref{Dircond} and $u_0$
be the solution to the same boundary value problem with $D$ replaced by the
empty set, and let $v$ and $v_0$ be the solutions to the same boundary value
problems with $\psi$ replaced by $\varphi$. Then
$$\langle(F(D)-F_0)\psi,\varphi\rangle=\int_{\Omega}(1+(a-1)\chi_D)\nabla u\cdot\nabla v
-\int_{\Omega}\nabla u_0\cdot\nabla v_0,$$
which is clearly symmetric, thus the self-adjointness follows.

By self-adjointness,
it is enough to show that there exist positive
constants $C_2$ and $\alpha_2$, depending on $N$, $m$, $\beta$ and $a$ only, so that
\begin{equation}\label{claim}
\left\|(F(D)-F_0)\frac{f_{jp}}{\sqrt{1+j}}\right\|_{H^{-1/2}(\partial\Omega)}\leq
C_2\exp(-\alpha_2(1+j)),
\end{equation}
for any $j$ and $p$. In fact, by \eqref{claim}, we infer that
$$\begin{array}{rl}
\left|\left\langle(F(D)-F_0)\frac{f_{jp}}{\sqrt{1+j}},\frac{f_{kq}}{\sqrt{1+k}}\right\rangle\right|&\leq
\left\|(F(D)-F_0)\frac{f_{jp}}{\sqrt{1+j}}\right\|_{H^{-1/2}(\partial\Omega)}\\&\leq
C_2\exp(-\alpha_2(1+j)),\end{array}$$
where we have used the fact that $\left\|\frac{f_{kq}}{\sqrt{1+k}}\right\|_{H^{1/2}(\partial\Omega)}=1$.
Since $F(D)-F_0$ is self-adjoint, we can reverse the role of $j$ and $k$ and thus obtain
$$\left|\left\langle(F(D)-F_0)\frac{f_{jp}}{\sqrt{1+j}},\frac{f_{kq}}{\sqrt{1+k}}\right\rangle\right|\leq
C_2\exp(-\alpha_2(1+k))$$
as well, and so \eqref{mainestimate} immediately follows from \eqref{claim}.

Let $u_{jp}(D)$ be the solution to \eqref{incleq}-\eqref{Dircond} with the
boundary datum $\psi$ replaced by $f_{jp}$, let $u_{jp}$ be defined as in
\eqref{polyharm} and let $v_{jp}=u_{jp}(D)-u_{jp}$.
Since every $D\in X$ is compactly contained in $B_N(0,4/5)$, we can find a
constant $C_3$, which depends only on $N$, so that for any $D\in X$
$$\left\|(F(D)-F_0)\frac{f_{jp}}{\sqrt{1+j}}\right\|_{H^{-1/2}(\partial\Omega)}\leq
C_3\left(\int_{\Omega\backslash\overline{B_N(0,4/5)}}
\|\nabla v_{jp}\|^2\right)^{1/2}.$$
In fact, there exists a constant $C_3$, depending on $N$ and $4/5$ only,
so that for every $\psi\in H^{1/2}(\partial\Omega)$ there exists $\tilde{\psi}\in H^1(\Omega)$
with the properties that $\tilde{\psi}|_{\partial\Omega}=\psi$, $\tilde{\psi}\equiv 0$ on $B_N(0,4/5)$ and
$$\|\tilde{\psi}\|_{H^1(\Omega)}=
\|\tilde{\psi}\|_{H^1(\Omega\backslash\overline{B_N(0,4/5)})}\leq C_3\|\psi\|_{H^{1/2}(\partial\Omega)}.$$

We have that $v_{jp}$ solves the boundary value problem
\begin{equation}\label{eq1}
\left\{\begin{array}{ll}
\mathrm{div}((1+(a-1)\chi_D)\nabla v_{jp})=-\mathrm{div}((a-1)\chi_D\nabla u_{jp}) 
&\text{in }\Omega,\\
v_{jp}=0&\text{on }\partial\Omega.
\end{array}\right.
\end{equation}

By the weak formulation of \eqref{eq1}
and the fact that $\overline{D}\subset B_N(0,4/5)$,
we have that for every $\psi\in H^{1/2}(\partial\Omega)$
$$\langle(F(D)-F_0)f_{jp},\psi\rangle=\int_{\Omega}\nabla v_{jp}\cdot\nabla\tilde{\psi},$$
thus
$$\begin{array}{rl}
|\langle(F(D)-F_0)f_{jp},\psi\rangle|&\leq
\left(\int_{\Omega\backslash\overline{B_N(0,4/5)}}
\|\nabla v_{jp}\|^2\right)^{1/2}\!\!
\left(\int_{\Omega\backslash\overline{B_N(0,4/5)}}
\|\nabla \tilde{\psi}\|^2\right)^{1/2}\\&\leq
C_3\left(\int_{\Omega\backslash\overline{B_N(0,4/5)}}
\|\nabla v_{jp}\|^2\right)^{1/2}
\|\psi\|_{H^{1/2}(\partial\Omega)}.\end{array}$$

Again by the weak formulation of \eqref{eq1}, we have that
$$\begin{array}{rl}\int_{\Omega}
\|\nabla v_{jp}\|^2&\leq \max\{1,1/a\}
\int_{\Omega}(1+(a-1)\chi_D)\|\nabla v_{jp}\|^2\\&=
\max\{1,1/a\}(-\int_D(a-1)\nabla u_{jp}\cdot\nabla v_{jp})\\&\leq
\max\{1,1/a\}|a-1|\left(\int_{\Omega}
\|\nabla v_{jp}\|^2\right)^{1/2}\left(\int_{D}
\|\nabla u_{jp}\|^2\right)^{1/2}.\end{array}$$

From here, using again the fact that $D$ is contained in $B_N(0,4/5)$,
it is easy to infer that there exists a constant $C_4$ depending on $a$
only so that
$$\left(\int_{\Omega\backslash\overline{B_N(0,4/5)}}
\|\nabla v_{jp}\|^2\right)^{1/2}\leq C_4 \left(\int_{B_N(0,4/5)}
\|\nabla u_{jp}\|^2\right)^{1/2}$$
and so \eqref{claim} is proved by using \eqref{expdecayball}.\cvd

\subsubsection{Neumann-to-Dirichlet case}

\proof{ of Proposition~\ref{inclprop} (Neumann-to-Dirichlet case).}
For what concerns the Neumann-to-Dirichlet case, the proposition can be proved
as an application of the abstract theorem (see for a similar procedure the
proof of
the ``insulating crack \& Neumann-to-Dirichlet case'' for the inverse crack problem,
page~\pageref{insNDcase}).
However, it is also possible to argue
in the following way. We observe that, for every $D\in X$,
$\mathcal{N}(D)$ and $\tilde{\Lambda}(D)$, the restriction of $\Lambda(D)$ to
${}_0H^{1/2}(\partial\Omega)$, are inverse to each other. Since we have already established
the Dirichlet-to-Neumann case, there exist $\varepsilon_1>0$ and $C$, depending
on $N$, $M$, $\beta$ and $a$ only, so that for any $\varepsilon$, $0<\varepsilon<\varepsilon_1$,
there exists $D_1$ and $D_2$ in $X$ so that \eqref{firstpart} and \eqref{DNpart} are satisfied.
Using the identity
\begin{equation}\label{Aleid}
\mathcal{N}(D_1)-\mathcal{N}(D_2)=\mathcal{N}(D_2)(\tilde{\Lambda}(D_2)-\tilde{\Lambda}(D_1))\mathcal{N}(D_1),
\end{equation}
we infer that
\begin{equation}\label{interpest}
\|\mathcal{N}(D_1)-\mathcal{N}(D_2)\|\leq\|\mathcal{N}(D_2)\|\,\|\tilde{\Lambda}(D_2)-\tilde{\Lambda}(D_1)\|\,
\|\mathcal{N}(D_1)\|,
\end{equation}
where the natural norms have been used.
We remark that there exists a positive constant $C_5$, depending on
$N$, $M$, $\beta$ and $a$ only, so that for any $D\in X$ we have
\begin{equation}\label{unifest}
\|\mathcal{N}(D)\|\leq C_5.
\end{equation}
So \eqref{NDpart} is proved by \eqref{DNpart}, \eqref{interpest} and \eqref{unifest}.\cvd

\subsubsection{Experimental measurements case}

\proof{ of Proposition~\ref{inclpropexpmeas}.}
For the basic properties of the problem \eqref{expmeaspbm}, we shall always
refer to \cite{Som e Che e Isa}.
Let $H=H^1(\Omega)\times \mathbb{R}^L$.
For any $D\in X$, and for any $(u,U)$ and $(w,W)\in H$, let
$$B_{D}((u,U),(w,W))=\int_{\Omega}(1+(a-1)\chi_D)\nabla u\cdot\nabla w+
\sum_{l=1}^L\int_{e_l}(u-U_l)(w-W_l).$$
Then, for any $\eta\in{}_0H^{-1/2}(\partial\Omega)$ and any $I\in\mathbb{R}^L$
so that $\sum_{l=1}^LI_l=0$, we have that there exist, and it is unique up to an additive constant,
a couple $(u,U)\in H$ so that we have
\begin{equation}\label{eqform}
B_{D}((u,U),(w,W))=\eta(w|_{\partial\Omega})+\sum_{l=1}^LI_lW_l,\quad\text{for any }(w,W)\in H.
\end{equation}
Furthermore, we have that 
there exists a constant $C_6$, depending on $N$, $a$, $Z$ and the electrodes only,
so that
\begin{equation}\label{coerc}
\|\nabla u\|_{L^2(\Omega)}\leq C_6(\|\eta\|_{{}_0H^{-1/2}(\partial\Omega)}+\|I\|).
\end{equation}

We also remark that if $(u,U)$ solves \eqref{eqform}, then
$$\mathrm{div}((1+(a-1)\chi_D)\nabla u)=0\quad\text{in }\Omega.$$
Therefore, we can rewrite the equation \eqref{eqform} on the boundary as follows. We denote
$\phi=\frac{\partial u}{\partial\nu}|_{\partial\Omega}\in{}_0H^{-1/2}(\partial\Omega)$
and we recall that
$\mathcal{N}(D)$ is the Neumann-to-Dirichlet map associated to the inclusion $D$.
Then, we deduce by straightforward computations that $(u,U)$ satisfies \eqref{eqform} if and only if 
\begin{equation}\label{Uexp}
U_l=\frac{z_l}{|e_l|}I_l+\frac{1}{|e_l|}\int_{e_l}u,
\quad\text{for any }l=1,\ldots,L,
\end{equation}
and the following equation holds in ${}_0H^{-1/2}(\partial\Omega)$
\begin{equation}\label{boundaryeqform}
\phi+\sum_{l=1}^L\frac{1}{z_l}\left(\mathcal{N}(D)\phi-
\frac{1}{|e_l|}\int_{e_l}\mathcal{N}(D)\phi\right)\chi_{e_l}
=\eta+\sum_{l=1}^L\left(\frac{I_l}{|e_l|}\chi_{e_l}\right).
\end{equation}

Let $\mathcal{K}(D):{}_0H^{-1/2}(\partial\Omega)\mapsto{}_0H^{-1/2}(\partial\Omega)$ be the
operator defined as follows. For any $\phi\in {}_0H^{-1/2}(\partial\Omega)$,
$$\mathcal{K}(D)\phi=\sum_{l=1}^L\frac{1}{z_l}\left(\mathcal{N}(D)\phi-
\frac{1}{|e_l|}\int_{e_l}\mathcal{N}(D)\phi\right)\chi_{e_l}.$$

We have that $\mathcal{K}(D)$ is a compact linear operator. Since the equation
\eqref{eqform} admits, up to additive constants, a unique solution, we can infer
that \eqref{boundaryeqform} is uniquely solvable, therefore the
operator $Id+\mathcal{K}(D)$ is invertible, where $Id$ denotes the identity operator.
Using \eqref{coerc}, this inverse,
which we shall denote with $\tilde{\mathcal{K}}(D)$,
satisfies
\begin{equation}\label{unifest2}
\|\tilde{\mathcal{K}}(D)\|=\|(Id+\mathcal{K}(D))^{-1}\|\leq C_7,
\end{equation}
where the constant $C_7$ depends
on $N$, $a$, $Z$ and the electrodes only.

For any given current pattern $I$, that is $I\in \mathbb{R}^L$ so that
$\sum_{l=1}^LI_l=0$, we can define $\tilde{I}=\sum_{l=1}^L\left(\frac{I_l}{|e_l|}\chi_{e_l}\right)\in
{}_0H^{-1/2}(\partial\Omega)$. Furthermore, there exists a constant $C_8$, depending on $N$ and
the electrodes only, so that, for any $I\in \mathbb{R}^L$ satisfying
$\sum_{l=1}^LI_l=0$, we have
\begin{equation}\label{ItildeI}
\|\tilde{I}\|_{{}_0H^{-1/2}(\partial\Omega)}\leq C_8\|I\|.
\end{equation}
As it is shown in \cite{Som e Che e Isa},
we have that $u$ solves our direct problem \eqref{expmeaspbm}
for a given current pattern $I$
if and only if
\eqref{eqform} is satisfied with $\eta=0$. Therefore, if
we take $I\in\mathbb{R}^L$ so that
$\sum_{l=1}^LI_l=0$, we have that $R(D)I=V$ where, for any $l=1,\ldots,L$,  
\begin{equation}\label{R-Kconn}
V_l=\int_{e_l}\mathcal{N}(D)(\tilde{\mathcal{K}}(D)\tilde{I})+c|e_l|,
\end{equation}
where $c$ is a constant which can be computed by imposing the condition
that $\sum_{l=1}^LV_l=0$, that is
\begin{equation}\label{constantc}
c=-\frac{\sum_{l=1}^L\int_{e_l}\mathcal{N}(D)(\tilde{\mathcal{K}}(D)\tilde{I})}{\sum_{l=1}^L|e_l|}.
\end{equation}

In order to establish Proposition~\ref{inclpropexpmeas},
we observe that, by Proposition~\ref{inclprop}, we can find constants $\varepsilon_1>0$ and $C$, which depend
on $N$, $m$, $\beta$ and $a$ only, so that for any
$\varepsilon$, $0<\varepsilon<\varepsilon_1$,
there exists $D_1$ and $D_2$ in $X$ so that \eqref{firstpart} and \eqref{NDpart} are satisfied.
We show that these inclusions $D_1$ and $D_2$ provide us with the instability example also in
the experimental measurements case.
Let us evaluate
the norm of $R(D_1)-R(D_2)$. Therefore, we take $I\in\mathbb{R}^L$ so that
$\sum_{l=1}^LI_l=0$ and we evaluate
$\|(R(D_1)-R(D_2))I\|$. We recall that we have posed $R(D_1)[1]=R(D_2)[1]=0$.
By \eqref{R-Kconn} and \eqref{constantc},
we have that
$$\|(R(D_1)-R(D_2))I\|\leq C_9\|\mathcal{N}(D_1)(\tilde{\mathcal{K}}(D_1)\tilde{I})-
\mathcal{N}(D_2)(\tilde{\mathcal{K}}(D_2)\tilde{I})\|_{L^2(\partial\Omega)},$$
where $C_9$ depends on $N$ and the electrodes only. Thus, we can find a constant $C_{10}$,
depending on $N$ and $C_{9}$ only, so that
$$\begin{array}{rcl}
\|(R(D_1)-R(D_2))I\| &\leq &C_{10}\|\mathcal{N}(D_1)(\tilde{\mathcal{K}}(D_1)-\tilde{\mathcal{K}}(D_2))
\tilde{I}\|_{H^{1/2}(\partial\Omega)}+\\
&& C_{10}\|(\mathcal{N}(D_2)-\mathcal{N}(D_1))(\tilde{\mathcal{K}}(D_2)
\tilde{I})\|_{H^{1/2}(\partial\Omega)},
\end{array}$$
and, by \eqref{unifest}, \eqref{unifest2} and \eqref{ItildeI}, we also deduce that
\begin{equation}\label{Rest}
\begin{array}{rcl}
\|(R(D_1)-R(D_2))I\|&\leq&
C_5C_8C_{10}\|\tilde{\mathcal{K}}(D_1)-\tilde{\mathcal{K}}(D_2)\|\,
\|I\|+\\
& & C_7C_8C_{10}\|\mathcal{N}(D_2)-\mathcal{N}(D_1)\|\,\|I\|.
\end{array}
\end{equation}

It remains to evaluate the term
$\|\tilde{\mathcal{K}}(D_1)-\tilde{\mathcal{K}}(D_2)\|$. We proceed as follows.
Using an identity analogous to \eqref{Aleid} applied to the operators $Id+\mathcal{K}$
and $\tilde{\mathcal{K}}$, and recalling \eqref{unifest2}, we obtain that \begin{equation}\label{K1-K2}
\|\tilde{\mathcal{K}}(D_1)-\tilde{\mathcal{K}}(D_2)\|\leq C_7^2
\|\mathcal{K}(D_1)-\mathcal{K}(D_2)\|\leq C_{11}\|\mathcal{N}(D_1)-\mathcal{N}(D_2)\|,
\end{equation}
where $C_{11}$ depends on $N$, $a$, $Z$ and the electrodes only.

And so the conclusion immediately follows by coupling \eqref{Rest} with \eqref{K1-K2} and using
\eqref{NDpart}.\cvd

\subsection{Inverse crack problem}\label{cracksubproof}

The proof of Proposition~\ref{invcrackprop}
follows directly from the abstract theorem stated in
Theorem~\ref{abstractthm}. We just need to check that all the hypotheses of
Theorem~\ref{abstractthm} are satisfied. Therefore, the proof
is divided into two steps, each corresponding to one of the hypotheses of
Theorem~\ref{abstractthm}.

\proof{ of Proposition~\ref{invcrackprop} - First step.}
First, let $x_0\in X$ be $\overline{B'_{N-1}(0,1/2)}$. Then,
by Proposition~\ref{discretesetprop},
$X$ satisfies assumption~\textit{\ref{Xassumpt}}\textnormal{)}
of Theorem~\textnormal{\ref{abstractthm}}, with constants $\varepsilon_0$ and
$C_1$
depending
on $N$, $m$ and $\beta$ only, and constant $\alpha_1=(N-1)/m$.
We recall also that $\sigma\subset B_N(0,4/5)$ for any $\sigma\in X$.\cvd

For what concerns the second step,
we turn our attention to assumption~\textit{\ref{HFassumpt}}\textnormal{)}
of Theorem~\textnormal{\ref{abstractthm}}. Each case, corresponding to
operators $\Lambda_i$ and $\mathcal{N}_i$, $i=1,2$, should be treated separately. We limit ourselves to
two cases, namely the cases corresponding to $\mathcal{N}_1$ and $\Lambda_2$,
in order to show the main points of the proof,
and we leave the details concerning the other two cases to the reader.

\subsubsection{Insulating crack \& Neumann-to-Dirichlet case}\label{insNDcase}

\proof{ of Proposition~\ref{invcrackprop} - Second step (Insulating crack \& Neu\-mann-to-Diri\-chlet case).}
First, we notice that $u$ is a solution to \eqref{inscrack}-\eqref{currentbdrycond}
if and only if $u\in H^1(\Omega\backslash\sigma)$, $\int_{\partial\Omega}u=0$, and
$$\int_{\Omega\backslash\sigma}\nabla u\cdot\nabla w=\eta(w|_{\partial\Omega}),\quad
\text{for any }w\in H^1(\Omega\backslash\sigma).$$

We observe that for any 
$\sigma\in X$, $\mathcal{N}_1(\sigma)$ is a bounded and linear operator between
${}_0H^{-1/2}(\partial\Omega)$ and its dual. Hence we take
$H$ to be ${}_0H^{-1/2}(\partial\Omega)$
and $F:X\mapsto\mathcal{L}(H,H')$ to be defined as $F(\sigma)=\mathcal{N}_1(\sigma)$ for any
$\sigma\in X$.
With $F_0$ we denote in an analogous way the Neumann-to-Dirichlet
map related to \eqref{inscrack}-\eqref{currentbdrycond} 
when $\sigma=\emptyset$, that is the Neumann-to-Dirichlet map associated to
the body where no crack is present.
For any $\eta\in{}_0H^{-1/2}(\partial\Omega)\backslash\{0\}$, we define
\begin{equation}\label{defgammaNeu}
\gamma(\eta)=
\frac{\|\eta\|^2_{L^2(\partial\Omega)}}{\|\eta\|^2_{H^{-1/2}(\partial\Omega)}}.
\end{equation}

Referring to Proposition~\ref{PA1prop}, Example~\ref{ballex} and
\eqref{spherarm}, $\{v_k\}_{k\in\mathbb{N}}$,
the orthonormal basis of $H$ we shall employ, is given by
\begin{equation}\label{basisdefNeu}
\left\{\sqrt{j}f_{jp}:\ j\geq 1\text{ and }1\leq p \leq p_j\right\}
\end{equation}
with the natural order. We have that $\gamma(\sqrt{j}f_{jp})=j$,
for any $j$ and $p$.
Again by our remarks in Example~\ref{ballex}, we deduce that
$\#\{k\in\mathbb{N}:\ \gamma(v_k)\leq n\}\leq 2(1+n)^{N-1}$,
for any $n\in\mathbb{N}$,

For what concerns
\eqref{mainestimate}, we argue in this way.
We need a kind of self-adjointness of $F(\sigma)-F_0$ for every $\sigma\in X$. We have that
$$\langle(F(\sigma)-F_0)\eta,\phi\rangle=\langle(F(\sigma)-F_0)\phi,\eta\rangle$$
for any $\eta$, $\phi\in{}_0H^{-1/2}(\partial\Omega)$, where $\langle\cdot,\cdot\rangle$
is again the
duality pairing between $H'$ and $H$.
In fact, if $u$ solves
\eqref{inscrack}-\eqref{currentbdrycond}, $u_0$ solves
the same boundary value problem with $\sigma$ replaced by the
empty set, $v$ and $v_0$ solves the same boundary value
problems with $\eta$ replaced by $\phi$, then
$$\langle(F(\sigma)-F_0)\eta,\phi\rangle=
\int_{\Omega\backslash\sigma}\nabla v\cdot\nabla u
-\int_{\Omega}\nabla v_0\cdot\nabla u_0.$$

By the self-adjointness of the operator $F(\sigma)-F_0$, for any $\sigma\in X$,
in order to prove \eqref{mainestimate} we have to show that there 
exist positive constants $C_2$ and $\alpha_2$, which depend on $N$, $m$ and $\beta$ only,
so that, for any $j$ and $p$,
\begin{equation}\label{claim3}
\left\|(F(\sigma)-F_0)\sqrt{j}f_{jp}\right\|_{H^{1/2}(\partial\Omega)}\leq
C_2\exp(-\alpha_2j).
\end{equation}

We can find a constant $C_3$, depending on $N$ only, so that,
for any $\sigma\in X$,
$$\left\|(F(\sigma)-F_0)\sqrt{j}f_{jp}\right\|_{H^{1/2}(\partial\Omega)}\leq
C_3 \|v_{jp}\|_{H^1(\Omega\backslash\overline{B_N(0,4/5)})},$$
where $v_{jp}$ satisfies
\begin{equation}\label{auxiliarypbm}
\left\{\begin{array}{ll}
\Delta v_{jp}=0 &\text{in }\Omega\backslash\sigma,\\
\frac{\partial v_{jp}}{\partial\nu}=0&\text{on }\partial\Omega,\\
\frac{\partial v_{jp}}{\partial\nu}=-j^{-1/2}
\frac{\partial u_{jp}}{\partial\nu}&\text{on }
\partial\sigma,\\
\int_{\partial\Omega}v_{jp}=0,
\end{array}\right.
\end{equation}
with $u_{jp}$ given by formula \eqref{polyharm}. Since $\int_{\partial\Omega}v_{jp}=0$,
a Poincar\'e type inequality implies that there exists a constant $C_4$,
depending on $N$ only, so that, for any $\sigma\in X$, we have
$$\left\|(F(\sigma)-F_0)\sqrt{j}f_{jp}\right\|_{H^{1/2}(\partial\Omega)}\leq C_4
\left(\int_{\Omega\backslash\overline{B_N(0,4/5)}}\|\nabla
v_{jp}\|^2\right)^{1/2}.$$
We can estimate the right hand side of the last equation as follows. We fix a cut-off
function
$\chi$ so that $\chi\in C^{\infty}_0(B_N(0,5/6))$, $0\leq \chi\leq 1$, $\chi\equiv 1$ on
$B_N(0,4/5)$. Without loss of generality, we can assume that for every
$x\in \mathbb{R}^N$, $\|\nabla\chi(x)\|\leq C_5$, $C_5$ being a constant depending on
$N$ only.
Let us observe that
\eqref{auxiliarypbm} means that for every $w\in H^1(\Omega\backslash\sigma)$ we have
$$\int_{\Omega\backslash\sigma}\nabla v_{jp}\cdot\nabla w=-
\int_{\Omega\backslash\sigma}j^{-1/2}\nabla u_{jp}\cdot\nabla(\chi w).$$
Then, by taking $w=v_{jp}$, we infer that
$$\int_{\Omega\backslash\sigma}\|\nabla
v_{jp}\|^2=-\int_{\Omega\backslash\sigma}j^{-1/2}\nabla u_{jp}\cdot\nabla(\chi v_{jp}).$$
Straightforward computations allow us to prove that there exists a constant $C_6$,
depending on $N$ only, so that
$$\left(\int_{\Omega\backslash\sigma}\|\nabla
v_{jp}\|^2\right)^{1/2}\leq C_6\left(\int_{B_N(0,5/6)}\|\nabla
u_{jp}\|^2\right)^{1/2}.$$
Then we can conclude using \eqref{expdecayball}.\cvd

\subsubsection{Conducting crack \& Dirichlet-to-Neumann case}

\proof{ of Proposition~\ref{invcrackprop} - Second step (Conducting crack \& Diri\-chlet-to-Neumann
case).}
We begin with a description of the weak formulation of the boundary value problem
\eqref{condcrack}-\eqref{voltbdrycond}. With $H^1_{const}(\Omega,\sigma)$
we denote the subspace of $H^1(\Omega)$ functions which are constant on $\sigma$.
For any $c\in\mathbb{R}$, we set $H^1_c(\Omega,\sigma)$ as the subset of
$H^1(\Omega)$ functions which are equal to the constant $c$ on $\sigma$.
For any $c\in\mathbb{R}$, we have that there exists and it is unique a solution to the
following boundary value problem
\begin{equation}\label{cpbm}
\left\{\begin{array}{ll}
\Delta u_c=0 &\text{in }\Omega\backslash\sigma,\\
u_c=c &\text{on }\sigma,\\
u_c=\psi&\text{on }\partial\Omega,
\end{array}\right.
\end{equation}
that is a function $u_c\in H^1_c(\Omega,\sigma)$ so that $u_c|_{\partial\Omega}=\psi$
and that
$$\int_{\Omega\backslash\sigma}\nabla u_c\cdot\nabla w=0,\quad\text{for
any }w\in H^1_0(\Omega)\cap H^1_0(\Omega,\sigma).$$
Given $u_c$, solution to \eqref{cpbm}, we can define
$\frac{\partial u_c}{\partial\nu}|_{\partial\Omega}\in H^{-1/2}(\partial\Omega)$
as follows
$$\langle\frac{\partial u_c}{\partial\nu}|_{\partial\Omega},\varphi\rangle
=\int_{\Omega\backslash\sigma}\nabla u_c\cdot\nabla\tilde{\varphi},$$
where $\varphi\in H^{1/2}(\partial\Omega)$ and $\tilde{\varphi}$ is any
$H^1_0(\Omega,\sigma)$ function so that
$\tilde{\varphi}|_{\partial\Omega}=\varphi$.

We claim that there exists a unique $c\in\mathbb{R}$ so that
$\langle\frac{\partial u_c}{\partial\nu}|_{\partial\Omega},1\rangle=0$, that is
existence and uniqueness of a solution to
\eqref{condcrack}-\eqref{voltbdrycond}.

We have that $u$ solves \eqref{condcrack}-\eqref{voltbdrycond} if and only if
$u\in H^1_{const}(\Omega,\sigma)$ so that $u|_{\partial\Omega}=\psi$
and that
$$\int_{\Omega\backslash\sigma}\nabla u\cdot\nabla w=0,\quad\text{for
any }w\in H^1_0(\Omega)\cap H^1_{const}(\Omega,\sigma).$$
If we take $\tilde{\psi}$ to be any $H^1_0(\Omega,\sigma)$ function so that
$\tilde{\psi}|_{\partial\Omega}=\psi$, we have that $u$ solves
\eqref{condcrack}-\eqref{voltbdrycond} if and only if
$\tilde{u}=u-\tilde\psi$ belongs to $H^1_0(\Omega)\cap H^1_{const}(\Omega,\sigma)$
and satisfies
$$\int_{\Omega\backslash\sigma}\nabla\tilde{u}\cdot\nabla w=-
\int_{\Omega\backslash\sigma}\nabla\tilde{\psi}\cdot\nabla w,\quad\text{for
any }w\in H^1_0(\Omega)\cap H^1_{const}(\Omega,\sigma).$$
Standard elliptic theory provides us with existence and uniqueness of such a solution.
By the property
$\langle\frac{\partial u}{\partial \nu}|_{\partial\Omega},1\rangle=0$, we can infer that
$\frac{\partial u}{\partial\nu}|_{\partial\Omega}\in H^{-1/2}(\partial\Omega)$
can be also defined as
$$\langle\frac{\partial u}{\partial\nu}|_{\partial\Omega},\varphi\rangle
=\int_{\Omega\backslash\sigma}\nabla u\cdot\nabla\tilde{\varphi},$$
where $\varphi\in H^{1/2}(\partial\Omega)$ and $\tilde{\varphi}$ is any
$H^1_{const}(\Omega,\sigma)$ function so that
$\tilde{\varphi}|_{\partial\Omega}=\varphi$.

Now we can denote with $H$ the space $H^{1/2}(\partial\Omega)$, and we can
fix $\gamma$ as in
\eqref{defgammaDir} and the orthonormal basis as the one described in
\eqref{basisdefDir}.
The map $F:X\mapsto\mathcal{L}(H,H')$ is given by $F(\sigma)=\Lambda_2(\sigma)$, for
any $\sigma\in X$, and $F_0$ denotes the Dirichlet-to-Neumann map corresponding
to $\sigma=\emptyset$. We recall that the operator $F(\sigma)$ is self-adjoint for
any $\sigma\in X$, as well as $F_0$ is.

We proceed to verify \eqref{mainestimate} in this case.
First,
there exists a constant $C_7$, depending on $N$ only, so that,
for any $\sigma\in X$,
\begin{equation}\label{first}
\left\|(F(\sigma)-F_0)\frac{f_{jp}}{\sqrt{1+j}}\right\|_{H^{-1/2}(\partial\Omega)}\leq
C_7 \left(\int_{\Omega\backslash\overline{B_N(0,4/5)}}\|\nabla
v_{jp}\|^2\right)^{1/2} \end{equation}
where $v_{jp}=u_{jp}(\sigma)-\frac{u_{jp}}{\sqrt{1+j}}$, $u_{jp}(\sigma)$
being the solution to \eqref{condcrack}-\eqref{voltbdrycond}
with $\psi$ replaced by $\frac{f_{jp}}{\sqrt{1+j}}$ and $u_{jp}$
being as in \eqref{polyharm}.

Hence, $v_{jp}$ satisfies
\begin{equation}\label{auxiliarypbm2}
\left\{\begin{array}{ll}
\Delta v_{jp}=0 &\text{in }\Omega\backslash\sigma,\\
v_{jp}=0&\text{on }\partial\Omega,\\
v_{jp}=c-
\frac{u_{jp}}{\sqrt{1+j}}&\text{on }
\partial\sigma,\\
\langle\frac{\partial v_{jp}}{\partial \nu}|_{\partial\Omega},1\rangle=0,
\end{array}\right.
\end{equation}
where $c=u_{jp}(\sigma)|_{\sigma}$.
We notice that, if $\chi$ is the cut-off function previously defined in this subsection, then
$w_{jp}=(v_{jp}-c+\chi\frac{u_{jp}}{\sqrt{1+j}})\in H^1_0(\Omega,\sigma)$ and
$v_{jp}|_{\partial\Omega}=-c$. So,
$$\int_{\Omega\backslash\sigma}\nabla v_{jp}\cdot\nabla w_{jp}=\langle
\frac{\partial v_{jp}}{\partial \nu}|_{\partial\Omega},-c \rangle =0,$$
that is
$$\int_{\Omega\backslash\sigma} \nabla v_{jp}\cdot\nabla v_{jp}=
\int_{\Omega\backslash\sigma}\nabla v_{jp}\cdot\nabla (\chi\frac{u_{jp}}{\sqrt{1+j}}),$$
from which we easily deduce that
\begin{equation}\label{second}
\left(\int_{\Omega\backslash\sigma}\|\nabla
v_{jp}\|^2\right)^{1/2}\leq C_8\|u_{jp}\|_{H^1(B_N(0,5/6)},
\end{equation}
where $C_8$ depends on $N$ only.

So \eqref{mainestimate} is obtained by combining \eqref{first}, \eqref{second} and
\eqref{expdecayball} and the self-adjointness of the operator $F(\sigma)-F_0$.\cvd

\subsection{Inverse cavity problem, inverse surface crack problem
and inverse boundary material loss problem}\label{otherdefsubproof}

As we have already observed, the inverse problem of cavities can be treated in a way which is completely
analogous to the treatment of the inverse crack problem.

\proof{ of Proposition~\ref{surfcrackprop}.} It
can be obtained along the lines of the proof
of Proposition~\ref{invcrackprop}, with obvious modifications. In particular, the
reference point in $X$ is given by $\sigma_0$, the
orthonormal basis used are those described in Example~\ref{cutballex}, whereas the
reference operator is the one related to the domain $\Omega$, $\Omega$ as in
Example~\ref{cutballex}.\cvd

\proof{ of Proposition~\ref{bmlprop}.}
Also the arguments for the proof of Proposition~\ref{bmlprop} are simple modifications of
what we have used to prove Proposition~\ref{invcrackprop}, clearly making use of
the orthonormal basis described in Example~\ref{halfballex}.\cvd

\subsection{Inverse scattering problem}\label{scatsubproof}

The proof of Proposition~\ref{scatprop} is somehow different from the proofs of the analogous
propositions discussed previously. In fact, we can not prove
Proposition~\ref{scatprop} as a straightforward application of Theorem~\ref{abstractthm}.
Nevertheless, the procedure developed during the
proof of Theorem~\ref{abstractthm} can be adjusted in such a way to cover also the inverse
scattering case framework.
In the sequel, we limit ourselves to
the sound-soft case, the sound-hard case can be obtained with
minor adjustments. We shall point out the main differences between the sound-soft and the
sound-hard case and conclude the proof for the sound-hard case at the end of the subsection.
The proof of Proposition~\ref{scatprop} for the sound-soft case
will be divided in two steps.

\proof{ of Proposition~\ref{scatprop} - First step (Sound-soft case).}
First, we fix $x_0\in X$ to be equal to $\overline{B_{N}(0,1)}$ and we
observe that assumption~\textit{\ref{Xassumpt}}) of Theorem~\ref{abstractthm}
is satisfied, by Proposition~\ref{discretesetprop}, with constants $\varepsilon_0$ and
$C_1$, depending on $N$, $m$ and $\beta$ only, and constant
$\alpha_1=(N-1)/m$.\cvd

The second step deals with the main difference from the previous cases, which is as follows.
In the abstract theorem, we have a function $F$ which maps elements of a metric space
$X$ into elements of $\mathcal{L}(H,H')$, $H$ being a separable Hilbert space.
Now, fixed $a>0$, we define a map $F$ which associates to
each obstacle $D\in X$ a complex-valued function
defined on $S^{N-1}\times S^{N-1}$, namely $F(D)=\mathcal{A}_s(D)(\cdot,\cdot,a)$ or,
respectively,
$F(D)=\mathcal{A}_h(D)(\cdot,\cdot,a)$.
In the abstract theorem, fixed a suitable $F_0\in \mathcal{L}(H,H')$,
the operator $F(x)$, $x\in X$, was characterized
by the numbers $b_{k,l}=\langle(F(x)-F_0)v_k,v_l\rangle$, where $k$, $l\in\mathbb{N}$
and $\{v_k\}_{k\in\mathbb{N}}$ is a suitably chosen orthonormal basis of $H$. The
fundamental properties of such a characterization were summarized in
assumption~\textit{\ref{HFassumpt}}) of Theorem~\ref{abstractthm}. In particular,
the crucial property was a control on the asymptotic behaviour of the coefficients
$b_{k,l}$, which was provided by formulas \eqref{polygrowth} and \eqref{mainestimate}.
We shall obtain a completely analogous characterization by decomposing the
far-field pattern in spherical harmonics.

\proof{ of Proposition~\ref{scatprop} - Second step (Sound-soft case).}
We take $\{v_k\}_{k\in\mathbb{N}}$ as the orthonormal basis of $L^2(S^{N-1})$
described in Example~\ref{ballex}, precisely in \eqref{spherarm}, with the natural order.
Therefore, for each $k\in\mathbb{N}$, $v_k$ is a (real-valued) spherical
harmonic function on $S^{N-1}$. We set $\gamma(v_k)$ as the degree of the spherical
harmonic function $v_k$. We have that $\gamma(v_k)$ is an increasing sequence, with
respect to $k$, whose asymptotic behaviour satisfies \eqref{polygrowth} with constants
$C_2=2$ and $p=N-1$.

The decomposition of the far-field pattern in spherical harmonics is given by,
for any $(\hat{x},\omega,a)\in S^{N-1}\times S^{N-1}\times (0,\infty)$,
\begin{equation}\label{decomposition}
\mathcal{A}_s(D)(\hat{x},\omega,a)=\sum_{k,l}b_{k,l}(a)v_k(\hat{x})v_l(\omega),
\end{equation}
where the complex-valued coefficients $b_{k,l}(a)$ are given, for any $a\in (0,\infty)$,
by
\begin{equation}\label{decompcoeff}
b_{k,l}(a)=\int\!\!\!\int_{S^{N-1}\times S^{N-1}}\mathcal{A}_s(D)(\hat{x},\omega,a)
v_k(\hat{x})v_l(\omega)\mathrm{d}\hat{x}
\mathrm{d}\omega.
\end{equation}
Furthermore, we use the following characterization
\begin{equation}\label{decompcoeff2}
b_{k,l}(a)=\int_{S^{N-1}}\tilde{b}_k(\omega,a)
v_l(\omega)\mathrm{d}\omega,
\end{equation}
where
the complex-valued coefficients $\tilde{b}_k(\omega,a)$ are, for any
$\omega\in S^{N-1}$ and any
$a\in (0,\infty)$, the Fourier coefficients, with respect to the orthonormal basis
$\{v_k\}_{k\in\mathbb{N}}$,
of the far-field pattern $u^s_{\infty}(\cdot;\omega;a)$ corresponding to
the scattered field of the solution to \eqref{Helmeq}-\eqref{softbc}-\eqref{Somm}, that is
\begin{equation}\label{decompcoeff3}
\tilde{b}_k(\omega,a)=\int_{S^{N-1}}\mathcal{A}_s(D)(\hat{x},\omega,a)
v_k(\hat{x})\mathrm{d}\hat{x}.
\end{equation}

In the next lemma, we establish the asymptotic behaviour of the
coefficients $b_{k,l}$, which will play the role of the assumption stated in
\eqref{mainestimate}.

\begin{lem}\label{scatlemma}
Under the previous assumptions and definitions,
there exist positive constants $C_2$ and $\alpha_2$, depending
on $N$, $m$, $\beta$ and $I_N$ only, so that for any $D\in X$, for
any $a\in I_N$ and for any $(k,l)\in\mathbb{N}\times\mathbb{N}$,
we have
\begin{equation}\label{mainestscat}
|b_{k,l}(a)|\leq
C_2\exp(-\alpha_2\max\{\gamma(v_k),\gamma(v_l)\}),
\end{equation}
coefficient $b_{k,l}$ as in \eqref{decompcoeff}.
\end{lem}

\proof. First, we claim that there exist positive constants
$C_3$ and $\alpha_3$,  depending
on $N$, $m$, $\beta$ and $I_N$ only, so that for any $D\in X$, for
any $a\in I_N$, for any $\omega\in S^{N-1}$ and for any $k\in\mathbb{N}$, we have
\begin{equation}\label{mainclaim1}
|\tilde{b}_{k}(\omega,a)|\leq
C_3\exp(-\alpha_3\gamma(v_k)),
\end{equation}
$\tilde{b}_{k}$ defined by \eqref{decompcoeff3}.

By \eqref{mainclaim1} and \eqref{decompcoeff2}, we immediately
infer that, for any $k$, $l\in\mathbb{N}$,
\begin{equation}\label{mainclaim2}
|b_{k,l}(a)|\leq
|S^{N-1}|^{1/2}C_3\exp(-\alpha_3\gamma(v_k)),
\end{equation}
$|S^{N-1}|$ being the $(N-1)$-dimensional measure of $S^{N-1}$.

Then, we make use of the following \emph{reciprocity relation}, see for instance
\cite[Theorem~3.13]{Col e Kre98}. For any $D\in X$ and any $a\in(0,\infty)$ we
have
\begin{equation}\label{recrel}
\mathcal{A}_s(D)(\hat{x},\omega,a)=\mathcal{A}_s(D)(-\omega,-\hat{x},a),\quad\text{for
any }\hat{x},\omega\in S^{N-1}.
\end{equation} 
The reciprocity relation plays the role of self-adjointness for the elliptic operators
we have considered before and allows us, using \eqref{mainclaim2}, to easily
conclude the proof of the lemma. Therefore, what remains to be proven is the claim
in \eqref{mainclaim1}.

In order to prove \eqref{mainclaim1}, we begin with a uniform bound on the scattered
field. We notice that, for any $D\in X$, $D\subset B_N(9/5)$.
With a procedure which is analogous to the one used first 
in \cite[Lemma~2]{Isak92} and later in \cite{Ron02},
and using the fact that the scattered fields are radiating solutions to the Helmholtz
equation, we can find a constant $C_4$, depending on $N$, $m$, $\beta$ and
$I_N$ only,
so that, for any $D\in X$, any $\omega\in S^{N-1}$ and any $a\in I_N$, we have
\begin{equation}\label{inftybound}
|u^s(x;\omega;a)|\leq C_4\|x\|^{-(N-1)/2},\quad\text{for
any }x\in\mathbb{R}^N\backslash B_N(0,2),
\end{equation}
where $u^s$ is the scattered field corresponding to the solution to
\eqref{Helmeq}-\eqref{softbc}-\eqref{Somm}. We remark that only in the estimate above
the difference between the cases $N=2$ and $N=3$ shows up. We refer to \cite{Ron02} for
a detailed discussion about uniform estimates of decay at infinity for
radiating solutions to the Helmholtz equation.

By Theorem~2.14 in \cite{Col e Kre98}, we have that since $u^s$ is a radiating
solution to the Helmholtz equation, with coefficient $a>0$,
in $\mathbb{R}^N\backslash\overline{B_N(9/5)}$,
then, for any $x\in\mathbb{R}^N\backslash B_N(0,2)$,
\begin{equation}\label{decompscat}
u^s(x;\omega;a)=
\sum_k\hat{b}_k(\omega,a)H^{(1)}_{\gamma(v_k)}(\sqrt{a}\|x\|)v_k(x/\|x\|),
\end{equation}
where $\hat{b}_k$ are complex-valued coefficients given by
\begin{equation}\label{decompscatcoeff}
\hat{b}_kH^{(1)}_{\gamma(v_k)}(\sqrt{a}r)=\int_{S^{N-1}}u^s(r\hat{x};\omega;a)
v_k(\hat{x})\mathrm{d}\hat{x},\quad\text{for any }r\geq 2, \end{equation}
where, for any integer $n\geq 0$, $H^{(1)}_n$ denotes the \emph{Hankel function} of
first kind and order $n$.

Theorem~2.15 in \cite{Col e Kre98} provides us with the necessary link between
coefficients $\tilde{b}_k$ and $\hat{b}_k$. In fact, it holds that
\begin{equation}\label{link}
\tilde{b}_k(\omega,a)=(\pi/2)^{(N-3)/2}a^{-(N-1)/4}(-\rmi)^{\gamma_k+(N-1)/2}
\hat{b}_k(\omega,a).
\end{equation}

We combine \eqref{link} with \eqref{decompscatcoeff} and, by using \eqref{inftybound}, 
we obtain that there exists a constant $C_5$, depending on $C_4$ only, so that,
for any $k\in\mathbb{N}$, any $\omega\in S^{N-1}$ and any $a\in I_N$,
\begin{equation}\label{step1}
|\tilde{b}_k(\omega,a)|\leq C_5a^{-(N-1)/4}r^{-(N-1)/2}
|H^{(1)}_{\gamma(v_k)}(\sqrt{a}r)|^{-1},\quad\text{for any }r\geq 2.
\end{equation}
We choose $r$ as follows. For any $a$, $a\geq 1$, we take $r=2$, whereas for any
$a$, $0<a<1$, we pick $r=2/\sqrt{a}$. With this choice we infer that
for any $a\in I_N$, we have that
\begin{equation}\label{step2}
|\tilde{b}_k(\omega,a)|\leq C_6|H^{(1)}_{\gamma(v_k)}(\tilde{r})|^{-1}
\end{equation}
where $2\leq\tilde{r}\leq 2\max\{1,\overline{a}\}$ and the constant $C_6$ again
depends on $C_4$ only.
Then we can establish \eqref{mainclaim1} by using the well-known
asymptotic behaviour
of the Hankel functions. In fact, there exists a constant $C_7$, depending on $N$ and
$\overline{a}$ only, so that for any $\tilde{r}$, $2\leq\tilde{r}\leq 2\max\{1,\overline{a}\}$,
\begin{equation}\label{asymptHankel}
|H_n^{(1)}(\tilde{r})|^{-1}\leq \left\{\begin{array}{ll}
C_7 & \text{if }n=0,1,\\
C_7\left(\frac{\rme\tilde{r}}{2}\right)^n
(n-1)^{-(n-1)}&\text{for any }n\geq 2.
\end{array}\right.
\end{equation}
By a straightforward computation, \eqref{step2} together with \eqref{asymptHankel}
implies the validity of \eqref{mainclaim1}.\cvd

Now we have what is needed to prove Proposition~\ref{scatprop} in the sound-soft case. Let us just
notice that, for any $D\in X$ and any $a\in(0,\infty)$,
$$\|\mathcal{A}_s(D)(\cdot,\cdot,a)\|_{L^2(S^{N-1}\times S^{N-1})}=
\left(\sum_{k,l}|b_{k,l}(a)|^2\right)^{1/2}.$$
Then, with the same procedure used to prove Lemma~\ref{netlemma},
and keeping in mind the fact that we repeat the procedure $j$ times,
one for each $a_i$, $i=1,\ldots,j$,
we can find a constant $C_8$, depending on $N$, and the constants
$C_2$ and $\alpha_2$ of Lemma~\ref{scatlemma} only, so that, for any
$\delta$, $0<\delta<1/\rme$, there exists a subset $Y_{\delta}$
of $X$ with at most
$\exp(jC_8(-\log\delta)^{2N-1})$ elements so that for any $D\in X$ there exists
$\tilde{D}\in Y_{\delta}$ satisfying
$$\sup_{a\in \{a_1,\ldots,a_j\}}
\|(\mathcal{A}_s(D)-\mathcal{A}_s(\tilde{D}))(\cdot,\cdot,a)\|_{L^2(S^{N-1}\times S^{N-1})}
\leq\delta.$$
Then the conclusion of the proof of Proposition~\ref{scatprop} in the sound-soft case is immediate.\cvd

\proof{ of Proposition~\ref{scatprop} (Sound-hard case).}
We conclude this subsection sketching the proof for the sound-hard case.

First, for
any $\alpha$, $0<\alpha<1$, and $\beta>0$ we define
$Y_{(1,\alpha)\beta(1/2)}(S^{N-1}(0,1))$ as we have defined 
$Y_{m\beta(1/2)}(S^{N-1}(0,1))$, $m$ being an integer, with the only obvious modification
of replacing the $C^m$ norm with the $C^{1,\alpha}$ norm. Furthermore, we observe that
$Y_{(1,\alpha)\beta(1/2)}(S^{N-1}(0,1))$ satisfies assumption~\textit{\ref{Xassumpt}})
of Theorem~\ref{abstractthm} with constants $\varepsilon_0$ and
$C_1$, depending on $N$, $\alpha$ and $\beta$ only, and constant
$\alpha_1=(N-1)/(1+\alpha)$.

The difference between the sound-soft case and the sound-hard case relies in the
estimate contained in \eqref{inftybound}. With arguments which are analogous to the
ones used for the sound-soft obstacles, estimate \eqref{inftybound} can be proved
for sound-hard obstacles belonging to $Y_{m\beta(1/2)}(S^{N-1}(0,1))$, for any integer
$m\geq 2$, with a constant $C_4$ depending on $N$, $m$, $\beta$ and $I_N$ only,
and for sound-hard obstacles belonging to $Y_{(1,\alpha)\beta(1/2)}(S^{N-1}(0,1))$,
with $0<\alpha<1$, with a constant $C_4$ depending on $N$, $\alpha$,
$\beta$ and $I_N$ only.

Since the other part of the proof does not depend on the type of boundary conditions
used, the result follows, for sound-hard obstacles, for any $m\geq 2$. Some modifications
are needed to treat the case when $m=1$. We have that, for any
$0<\alpha<1$,
$Y_{(1,\alpha)\beta(1/2)}(S^{N-1}(0,1))\subset Y_{1\beta(1/2)}(S^{N-1}(0,1))$.
We apply the procedure described before to the set
$Y_{(1,\alpha)\beta(1/2)}(S^{N-1}(0,1))$ and, recalling Remark~\ref{remark1}, we
infer that for any $\gamma>0$ there exists a constant $\varepsilon_1>0$,
that depends on $N$, $j$, $\beta$, $I_N$, $\alpha$ and $\gamma$ only, so that
for any $\varepsilon$, $0<\varepsilon<\varepsilon_1$, we can find $D_3$ and $D_4$, both
belonging to $Y_{(1,\alpha)\beta(1/2)}(S^{N-1}(0,1))$, satisfying
$d(D_3,D_4)\geq \varepsilon$ and
$$
\sup_{a\in \{a_1,\ldots,a_j\}}
\|(\mathcal{A}_h(D_3)-\mathcal{A}_h(D_4))(\cdot,\cdot,a)\|_{L^2(S^{N-1}\times S^{N-1})}
\leq 2\exp(-\varepsilon^{-\alpha_1}),
$$
where $\alpha_1=\frac{N-1}{(1+\alpha)(2N-1+\gamma)}$.
We can choose, from the very beginning, $\alpha$ and $\gamma$ in such a way that
$\alpha$ and $\gamma$ depend on $N$ only and $(1+\alpha)(2N-1+\gamma)=2N$,
for instance we can take $\alpha=\frac{1}{4(2N-1)}$
and $\gamma=\frac{3(2N-1)}{4(2N-1)+1}$.
Thus the result is established also for
the case $m=1$ and sound-hard obstacles.\cvd

\end{document}